%% file: arxiv_version1.tex
\documentclass{article}
\usepackage[utf8]{inputenc}
\usepackage{amsmath,hyperref}
\usepackage{amssymb}
\usepackage{amsthm}
\usepackage{tikz-cd, circuitikz}
\usepackage{cleveref}
\usepackage[textsize=scriptsize]{todonotes}
\usepackage{xcolor}
\usepackage{tikz}
\usepackage{tikz-cd}
\usepackage{enumitem}
\usetikzlibrary{cd, calc, shapes, decorations.shapes, decorations.markings, decorations.pathreplacing, shapes.geometric, shapes.symbols,positioning}
\usepackage{pgfplots}
\usepackage{float}
\usepackage{url}
\usepackage{needspace}
\usepackage{graphicx}
\usepackage{subcaption}
\graphicspath{{./}{figs/}}
\newcommand{\subfigref}[2]{\figref{#1}(\subref{subfig:#1-#2})}

\newtheorem{theorem}{Theorem}[section]
\newtheorem{proposition}[theorem]{Proposition}
\newtheorem{lemma}[theorem]{Lemma}
\newtheorem{corollary}[theorem]{Corollary}
\newtheorem{conjecture}[theorem]{Conjecture}

\theoremstyle{definition}
\newtheorem{definition}[theorem]{Definition}

\newtheorem*{remark*}{Remark}

\newtheorem{example}[theorem]{Example}
\newtheorem{rems}[theorem]{Remarks}{\bfseries}{\itshape}

\AddToHook{env/lemma/begin}{\crefalias{theorem}{lemma}}
\AddToHook{env/definition/begin}{\crefalias{theorem}{definition}}
\AddToHook{env/corollary/begin}{\crefalias{theorem}{corollary}}
\AddToHook{env/proposition/begin}{\crefalias{theorem}{proposition}}

\newcommand{\lub}{\mathrm{l.u.b.}}

\newcommand{\RR}{\mathbb{R}}
\newcommand{\K}{\mathbb{K}}

\newcommand{\KK}{\mathcal{K}}

\definecolor{darkorchid}{rgb}{0.6,0.196,0.8}

\newcommand{\R}{\mathbb{R}}

\newcommand{\Z}{\mathbb{Z}}

\newcommand{\N}{\mathbb{N}}

\newcommand{\Nbar}{\overline{\N}}

\newcommand{\Rbar}{\overline{\R}}

\DeclareMathOperator{\lkmap}{lk}
\DeclareMathOperator{\simmap}{sim}
\DeclareMathOperator{\supp}{supp}
\DeclareMathOperator{\reach}{reach}
\DeclareMathOperator{\vertmap}{verts}
\DeclareMathOperator{\obname}{ob}
\DeclareMathOperator{\idx}{idx}

\DeclareMathOperator{\coords}{coords}
\DeclareMathOperator{\im}{im}

\DeclareMathOperator{\push}{push}
\renewcommand{\vec}[1]{\mathsf{Vec}(#1)}

\newcommand{\diagrams}{\mathcal{D}}
\newcommand{\cdiagrams}{\mathcal{CD}}
\newcommand{\cubes}[1]{\mathcal{#1}}
\newcommand{\clcubes}[1]{\overline{\cubes{#1}}}
\newcommand{\di}[1]{\dipaths{#1}{}{}}
\newcommand{\dipaths}[3]{\overrightarrow{P}_{#2}^{#3}(#1)}

\newcommand{\verts}[1]{{\vertmap({#1})}}

\newcommand{\thespace}[1]{| #1 |}
\newcommand{\vecx}{\mathbf{x}}
\newcommand{\vecy}{\mathbf{y}}
\newcommand{\vecw}{\mathbf{w}}

\newcommand{\vecp}{\mathbf{p}}
\newcommand{\vecq}{\mathbf{q}}
\newcommand{\vecz}{\mathbf{z}}
\newcommand{\veca}{\mathbf{a}}
\newcommand{\vecb}{\mathbf{b}}
\newcommand{\veczero}{\mathbf{0}}

\newcommand{\cl}[1]{\overline{#1}}

\newcommand{\propref}[1]{Property~\ref{prop:#1}.}

\newcommand{\figref}[1]{Figure~\ref{fig:#1}}
\newcommand{\lemref}[1]{Lemma~\ref{lem:#1}}

\newcommand{\exref}[1]{Example~\ref{ex:#1}}

\newcommand{\defref}[1]{Definition~\ref{def:#1}}

\newcommand{\stmtref}[1]{Statement~\emph{\ref{stmt:#1}}} 

\title{A Directed Cubical Framework for $k$-critical Bipersistence Modules}

\author{
  Robin Belton\footnote{Vassar College, Poughkeepsie, New York, United States}\\
  \texttt{rbelton@vassar.edu}
  \and Robyn Brooks \footnote{University of Utah, Utah, United States}\\
  \texttt{robyn.brooks@utah.edu}
   \and Lisbeth Fajstrup\footnote{Aalborg University, Aalborg, Denmark}\\
  \texttt{fajstrup@math.aau.dk}
  \and Brittany Terese Fasy\footnote{Montana State University, Bozeman, MT}\\
  \texttt{brittany.fasy@montana.edu}
  \and Elizabeth Vidaurre\footnote{Molloy University, Rockville Centre, NY}\\
  \texttt{evidaurre@molloy.edu}
}

\begin{document}

\maketitle


\begin{abstract}

In topological data analysis, multiparameter persistence modules ~(MPMs) are used to analyze data with multiple parameters.  Summarizing topological information for MPMs is a non-trivial task; one useful summary of bipersistence modules is the fibered barcode \cite{lesnick2015interactive,lesnick2025fast}. In this paper, we use a directed topological construction of the parameter space of a bifiltration to extend the fibered barcode.  Specifically, given a $k$-critical bifiltration $F$, we build an associated directed Euclidean cubical complex $C$ whose directed paths each correspond to a one-parameter subfiltration of $F$. The path space of $C$, $\di{C}$, is used to define the \emph{extended fibered barcode}, which contains all possible monotonically increasing paths in $C$. We show that there can be paths in $\di{C}$ which generate filtrations not obtainable along a line with positive slope, and that all lines with positive slope have a path in $\di{C}$ that generate an equivalent barcode. We introduce two types of equivalence for these restricted modules: combinatorial equivalence of diagrams, and simplicial equivalence of the filtrations generating the modules. We show that removing certain cubes from $C$ yields a subcomplex on which dihomotopic paths induce simplicially equivalent filtrations, and hence combinatorially equivalent diagrams.

\end{abstract}

\tableofcontents

\section{Introduction}
 Topological data analysis (TDA) integrates techniques from algebraic topology with data analysis to extract ``shape" from data.  A foundational method in TDA is \emph{persistent homology} (PH) which constructs a sequence of shapes from a data set, and measures how the homology classes evolve as we vary some parameter. The input is the data set and the output is a \emph{persistence barcode or diagram}. The persistence barcode is computable and stable\cite{cohensteiner2007stability,zomorodian2005computing}. PH has been successfully applied to various other domains including biological data~\cite{Belton2022Extremal,bhaskar2021topological,  Lawson2019, TorresEspin2021}, audio data~\cite{essl2025shiftaudio,SandersonComputational2017,  WojciechTopological24}, and the social sciences~\cite{Duchin2021persistence}. For more applications, see the database on non-theoretical uses of topology (DONUT)~\cite{DONUT}.

 However, there are many reasons that one might wish to analyze data along multiple parameters; for example, 1-d persistence can be quite sensitive to outliers. More generally, many data sets require analysis on more than one  parameter, making the use of multiparameter persistence (\cite{carlsson2007theory,frosini1999size}) particularly helpful.  Unfortunately, the generalization of a persistence barcode for modules over more than 1-dimension is not straightforward; although such barcodes can be defined, they fail to be complete invariants (see, for example, Section 4.1 of \cite{botnan2023introductionmultiparameterpersistence}). Thus, this leads us to look elsewhere for other partial invariants of multiparameter persistence modules.


  One such partial invariant is the \emph{fibered barcode}, defined for bipersistence modules in \cite{lesnick2015interactive,lesnick2025fast}.   To obtain the fibered barcode of a bipersistence module $M$, the parameter space $\R^2$ is fibered by lines $L$ with positive slope.  Each of these lines represents a restricted module $M_L$ which has a 1-d persistence barcode.  The fibered barcode for a module $M$ is defined as a map the set of all such lines to the space of persistence barcodes who's image contains the collection of all possible restrictions of $M$ obtained on any line with positive slope.  The fibered barcode is useful in that it is computable and provides a visual interpretation of data that varies along 2 parameters.  It also determines the \emph{matching distance} between two multiparameter persistence modules, which can be computed from the fibered barcode for bipersistence modules \cite{bapat2022computing,kerber2018exact}.

 In this paper, we look to generalize the fibered barcode for bipersistence modules by viewing the parameter space as a directed space, and looking at filtrations obtained along non-decreasing paths, i.e. dipaths.  Doing this allows for a larger class of possible obtained 1-d modules than can be obtained by restricting to lines with postive slope- it may be the case that there exists a restricted modules obtained along a dipath that is not obtainable on any line with positive slope.  Additionally, we look to leverage tools of directed topology, like directed homotopy of dipaths, to identify broad classes of barcodes included in the image of the fibered barcode.  We look to distinguish equivalence classes within the extended fibered barcode via two types of equivalence: combinatorial equivalence of persistence barcodes and simplicial equivalence of filtrations.
 
 There are several related ways that other mathematicians have investigated the fibered barcode or looked to understand multiparameter persistence modules by restriction to subposets. By viewing multiparameter persistence modules as sheafs, projected barcodes (\cite{BerkoukPetit2026,fernandes2025computation}) can be obtained by considering derived pushforwards of these modules along $\RR$.  Other restriction based invariants include \cite{DeyKimMemoli2024}, which develops a method for computing the generalized rank invariant of a bipersistence module by restricting the module to zigzag paths along the boundaries of intervals in the parameter space; this is generalized in \cite{DeyXin2024} to certain nice $d$-parameter persistence modules for $d>2$.  Invariants arising from general order-preserving embeddings of representation-finite posets are studied in \cite{AmiotBrustleHanson2026}. Bipath persistent homology, studying persistence modules generated by two filtrations with common endpoints, is introduced in \cite{AokiEscolarTada2025}; \cite{AlonsoLiu2025} provides a framework in which bipath persistence can be understood using zigzag persistence.  Finally, in \cite{sun2025path}, authors generalize the matching distance to monotone piecewise-linear paths in parameter space.  Our paper focuses on understanding the underlying filtrations and persistence modules obtained by restricting to any monotone path.

The outline of the paper is as follows.  In \Cref{sec:background}, we provide relevant definitions for directed spaces, Euclidean cubical complexes, simplicial filtrations, and persistence modules.  In \Cref{sec:dpsc}, we give a construction of a directed parameter space complex for a bipersistence module (\Cref{def:dpsc}), and describe the filtrations and modules obtained by restriction along paths in this complex.  \Cref{sec:moduleequivalence} introduces the two types of equivalence between restricted modules that we investigate: combinatorial equivalence of diagrams and simplicial equivalence of filtrations. We define the extended fibered barcode in \Cref{sec:fb} (\Cref{def:efb}) and obtain the result that, for $k$-critical bipersistence modules, the image of the fibered barcode is contained in the image of the extended fibered barcode (\Cref{cor:fb_in_efb}).  We then identify a subset of the parameter space complex in \Cref{sec:ec}, called the entrance complex (\Cref{def:ec}),  that has the property that dihomotopy equivalent paths within this complex will have simplicially equivalent barcodes (\Cref{lem:kcrit_dihomotopy}), and hence combinatorially equivalent diagrams.  Finally, we conclude with \Cref{sec:discussion}, a discussion that has conjectures and ideas of how to use this tool in the future.

\section{Background}\label{sec:background}
We set the scene by recalling definitions of directed spaces, Euclidean cubical complexes, and multiparameter persistence modules.

\subsection{Directed Topological Spaces}

We define directed spaces following the notation used in \cite{raussen2014homology,grandis2009diretced}.

\begin{definition}[Directed Topological Space]
Let $X$ be a topological space. A \emph{path} in $X$ is a continuous map from the unit interval \mbox{$\gamma:I=[0,1] \to X$}.  Let $X^I$ denote the space of all paths in $X$ equipped with the compact-open topology. If $\di{X}\subset X^I$ such that  
\begin{itemize}
    \item $\di{X}$ contains all constant paths,
    \item $\di{X}$ is closed under partial directed reparameterizations, where a partial directed reparameterization of a path $\gamma: I\to X$ is a composition $\gamma\circ f$, $f:I\to I$ such that $t\leq$ implies $f(t)\leq f(s)$, and
    \item $\di{X}$ is closed under path concatenation,
\end{itemize}

then the pair $(X,\di{X})$ is a \emph{directed space} or d-space. 
\end{definition}

For instance, the directed space with $X=\R^n$ and $\di{X}$ given by non-decreasing paths is called a \emph{Euclidean directed space}. Another space of interest is given by dipaths between two points. 

\begin{definition}[Space of Directed Paths from $\vecp$ to $\vecq$]
    If $\vecp, \vecq \in X$, the set of dipaths
from~$\vecp$ to $\vecq$, i.e., those paths~$\gamma\in\di{X}$ such that
$\gamma(0)=\vecp$ and $\gamma(1)=\vecq$ is denoted as $\dipaths{X}{\vecp}{\vecq}$.  This is a space with the compact open topology.

\end{definition}

Just as in the case with nondirected homotopy, pairs of paths can be dihomotopy equivalent.  In this paper, we assume that all dihomotopies between paths are path-homotopies, meaning that they are homotopies between paths with fixed endpoints.

\begin{definition}[Dihomotopy equivalence of dipaths]
  A \emph{dihomotopy} $\overrightarrow{H}:[0,1]\times [0,1] \to X$ between two dipaths $\gamma,\alpha$ in $\di{X}$ is a homotopy such that  
 \begin{itemize}
     \item $\gamma(t)=H(0,t)$ and $\alpha(t)=H(1,t)$, 
     \item $H(s,0)=\alpha(0)=\gamma(0)$ and $H(s,1)=\alpha(1)=\gamma(1)$ for all $s\in I$, and 
     \item for all $t\in [0,1]$, $H_t:=H(t,\cdot)$ is a dipath.
   \end{itemize} Two dipaths are \emph{dihomotopic} if their initial and terminal points are the same and there is a dihomotopy between them.  
\end{definition}

\subsection{Euclidean Cubical Complexes}

Let $(\R^d, \preceq)$ denote $\R^d$ under the product order.
Given two points $\vecp,\vecq \in \R^d$ with coordinates $p_i\in\R $ and $q_i\in \R$, respectively, and with $\vecp \preceq \vecq$, we denote
the \emph{closed and open hyperrectangle
from $\vecp$ to $\vecq$} as
\[
    [\vecp,\vecq] := \{\vecx \mid \vecp\preceq \vecx \preceq \vecq \},
    \text{ and }
    ]\vecp, \vecq[ := \{ \vecx \mid \vecp \prec \vecx \prec \vecq. \}
\]
We note here that $[\vecp,\vecq]$ is the \emph{closure} of $]\vecp, \vecq[$, and write:
$\cl{]\vecp, \vecq[}=[\vecp,\vecq]$.  The
dimension of~$[\vecp,\vecq]$ (and of~$]\vecp, \vecq[$)  is the number of nonzero
coordinates of the vector~$\vecp-\vecq$. A rectangle of dimension $n$ is an $n-$rectangle.
The zero-dimensional rectangle
is a single point, and we often denote the zero-rectangles ~$[\vecp,\vecp]$ and~$]\vecp,
\vecp[$ simply by~$\vecp$. 

\begin{definition}[Elementary Cube]
    
If $|\vecp_i-\vecq_i|\in\{0,1\}$ for all $i$, then we say that~$]\vecp,\vecq[$ is an \emph{unit
cube} or an \emph{elementary cube} (or simply, \emph{cube}, in this paper).
If we have two elementary cubes, $\tau$ and~$\sigma$, such that~${\tau} \subseteq \cl{\sigma}$, then we say that~$\tau$ is a face of~$\sigma$ and that~$\sigma$ is a coface of $\tau$.
\end{definition}

For any rectangular region in $\R^2$, we identify two subsets in the boundary of the rectangle as ``past" and ``future" components of boundary.


\begin{definition}[Past and Future Boundary]

The \emph{future boundary} of the closed rectangle $[\vecp,\vecq]\subset\R^2$, denoted $\partial^+([\vecp,\vecq])$, is defined as

\[\partial^+([\vecp,\vecq]):=[(p_1,q_2),\vecq]\cup[(q_1,p_2),\vecq];\]

The \emph{past boundary} of $[\vecp,\vecq]$, denoted $\partial^-([\vecp,\vecq])$, is defined as

\[\partial^-([\vecp,\vecq]):=[\vecp,(p_1,q_2)]\cup[\vecp,(q_1,p_2)].\]
    
\end{definition}

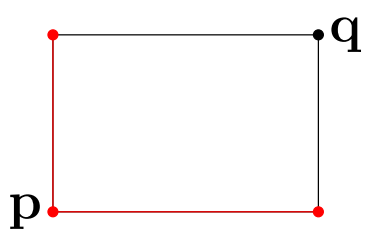
\begin{figure}[!ht]
\centering
\resizebox{1in}{!}{%
\begin{circuitikz}
\tikzstyle{every node}=[font=\fontsize{18.2pt}{23.7pt}\selectfont]
\draw  (6,13) rectangle (9,11);
\begin{scope}[transparency group, opacity=1]
\node at (6,13) [circ, color={rgb,255:red,255; green,0; blue,0}] {};
\end{scope}
\begin{scope}[transparency group, opacity=1]
\node at (9,13) [circ, color={rgb,255:red,0; green,0; blue,0}] {};
\end{scope}    \node [align = center, anchor = west] at (9,13) {%
 $\vecq$};
\begin{scope}[transparency group, opacity=1]
\node at (9,11) [circ, color={rgb,255:red,255; green,0; blue,0}] {};
\end{scope}
\begin{scope}[transparency group, opacity=1]
\node at (6,11) [circ, color={rgb,255:red,255; green,0; blue,0}] {};
\end{scope}   \node [align = center, anchor = east] at (6,11) {%
 $\vecp$};
\draw [ color={rgb,255:red,255; green,0; blue,0}, draw opacity=1, short] (6,11) -- (9,11);
\draw [ color={rgb,255:red,255; green,0; blue,0}, draw opacity=1, short] (6,11) -- (6,13);
\end{circuitikz}
}%
\label{future-past}\caption{The past boundary of the cube is in red and the future boundary is in black.}
\end{figure}

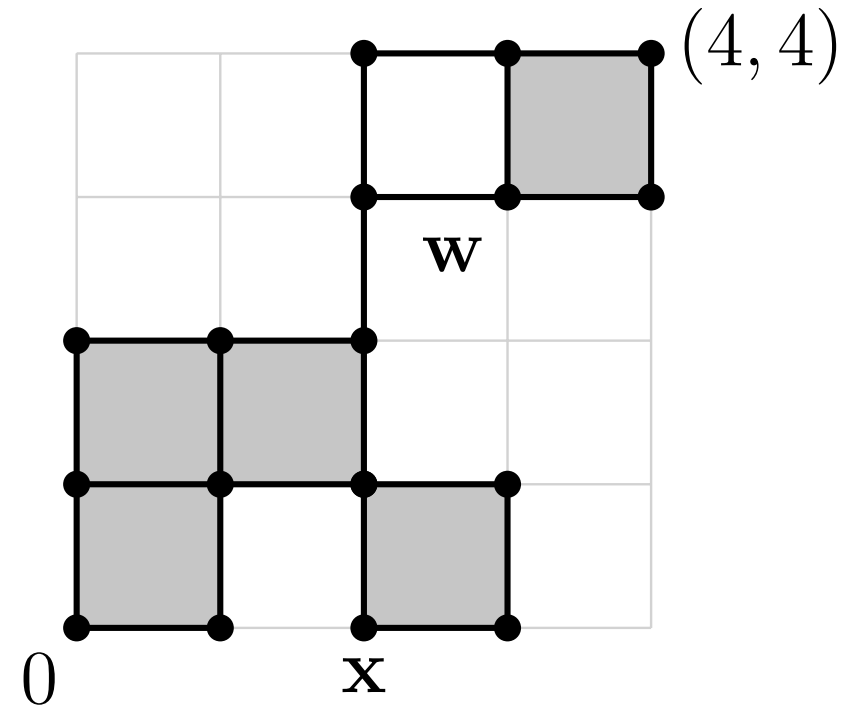
\begin{figure}[h!]
    \centering
    {
    \begin{tikzpicture}
        \input{figs/tixzfigures/ecc_example}  
    \end{tikzpicture}
     }
    \caption{Euclidean cubical complex in $\R^2$. The two-cubes are $[\veczero,(1,1)]$, $[(1,1),(2,2)]$,
    $[\vecx, (2,1)]$, and $[\vecw, (4,4)]$.
    }\label{fig:ecc}
\end{figure}

\begin{definition}[Euclidean Cubical Complex]
    A \emph{Euclidean cubical complex} (ECC) $(O,\mathcal{O})$ in $\R^d$ is a set $O\subset
\R^d$ along
with a finite
collection of (open) elementary cubes $\mathcal{O}$ in~$\R^d$ that is closed under the
face relation such that $\sqcup \mathcal{O}=O$; see
\figref{ecc}.
\end{definition}

Such cubical complexes
have been studied since the 1950s; see,
e.g.,~\cite{serre1951homologie,BrownOn81,jardine2002cubical,goubault1995geometry}.
We topologize this complex using the subspace topology, and use
$\thespace{\mathcal{O}}$
to denote $O$.\footnote{In general, given a set $\cubes{S}$ of cubes, we call
$\thespace{\cubes{S}}$ the `underlying
topological space.'}
For each $i \in \N$, let $\mathcal{O}_i$ denote the set of~$i$-cubes
in~$\mathcal{O}$. Often, we
denote the zero-cubes (vertices)  by
$\verts{\mathcal{O}}:=\mathcal{O}_0$.
Let~$\clcubes{O}$ denote the set of the closures of the open cubes in $\mathcal{O}$.

The d-spaces we consider are constructed from ECCs. Directed paths in a ECC are monotonically increasing paths from some initial vertex to a terminal vertex within the ECC, and we topologize the space of directed paths by endowing it with the compact-open topology. We refer to $(O, \mathcal{O}, \di{O})$ as a \emph{directed Euclidean cubical complex}.

\subsection{Filtrations}

In this section, we present definitions of filtrations and persistence modules.

\begin{definition}[Filtration]\label{def:filt}
    Let $(P, \leq)$ be a poset, and let $K$ be a simplicial complex.  We say that a
    function $F \colon P \to 2^K$ is a (monotone) \emph{filtration of $K$
    over~$P$} if the following
    conditions  are~met:
   
    \begin{enumerate}
        \item Subcomplex Condition: For all $x \in P$, $F(x)$ is a subcomplex of $K$.
        \item Monotonicity Condition: For each $x,x' \in P$ such that $x \leq x'$, $F(x)$ is a
            subcomplex of~$F(x')$.
        \item Full Complex Condition: $\bigcup_{x\in P}F(x) = K$.
    \end{enumerate}
    We refer to $P$ as the \emph{parameter space} for $F$.
\end{definition}

If $F$ is a filtration over $P$ and $P' \subset P$, then $F|_{P'}$ (the
restriction of $F$ to the poset $P'$) is a filtration of $\bigcup_{x\in P'}F(x)$.
When $P=[\vecp, \vecq ] \subseteq \R^d$, we say that the \emph{dimension} of the filtration is
the dimension of~$P$ (i.e., the number of nonzero entries in the vector
$\vecq-\vecp$).  We often refer to this type of filtration as
a~\emph{$k$-parameter filtration}, where $k$ is the dimension of $P$. \footnote{In this paper, we define many of our objects for $\R^d$, while prove results specifically for filtrations over $2$ dimensions.}

\begin{figure}[h!]
    \centering
    {\begin{tikzpicture}
        \input{figs/tixzfigures/bifiltration_grid}
    \end{tikzpicture}}
    \caption{A two-dimensional filtration. 
    }\label{fig:filtration}
\end{figure}
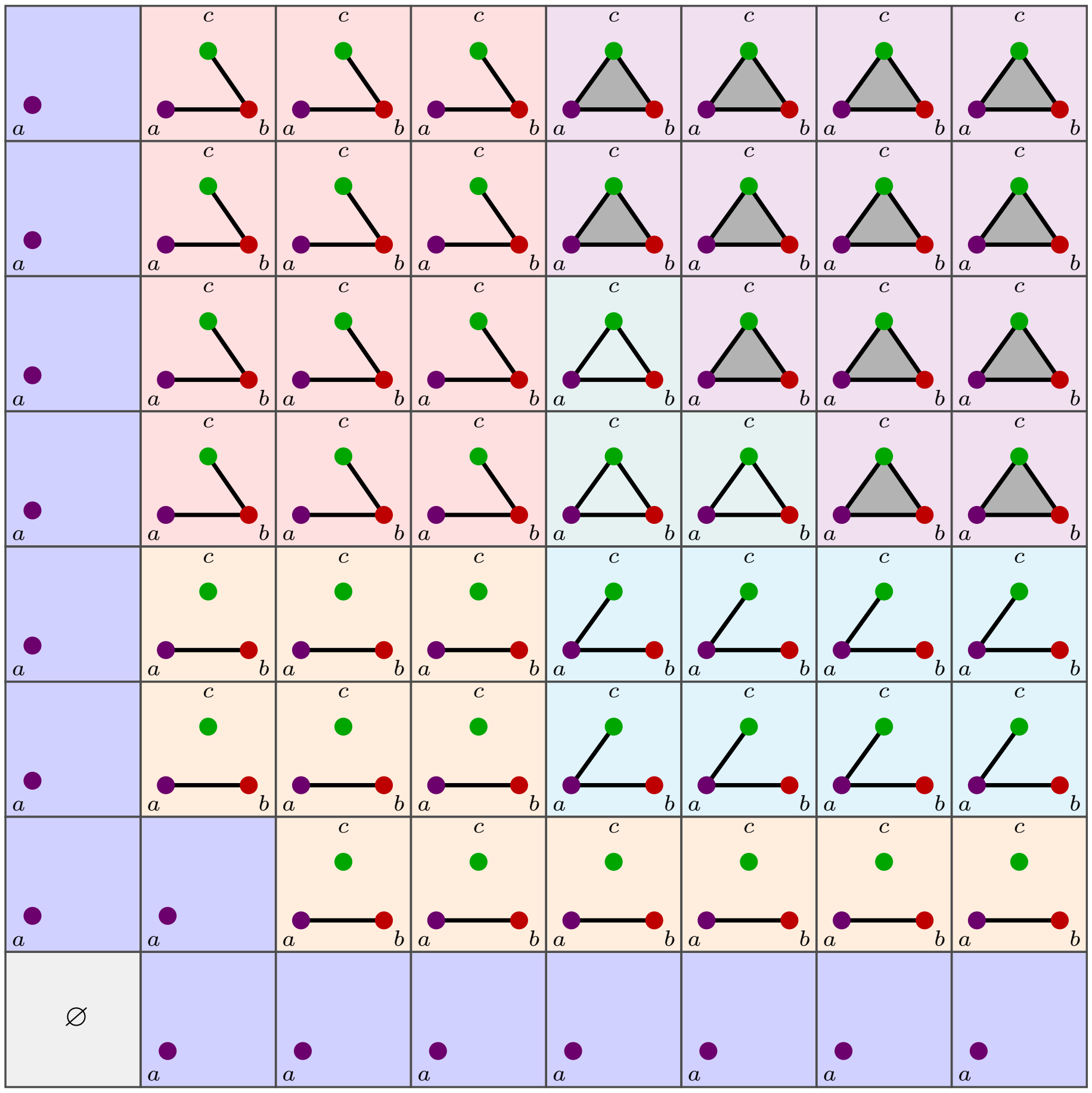

If $P$
has a unique terminal element, then, by the definition of filtration, the full
complex being filtered
is the complex at the terminal point of the~filtration. We formalize this as a
corollary to the definition:

\begin{corollary}[Filtrations of Posets with a Terminal
    Element]\label{cor:termposetfilt}
    If $P$ is a poset with a unique terminal element $\top$ and $F \colon P \to
    2^K$ is
    a filtration, then $F(\top)=K$.
\end{corollary}

The support of a simplex is the set of filtration parameters for which the
simplex is present:

\begin{definition}[Entrance coordinate, Support of a Simplex]
    Given a filtration $F$ of a simplicial complex $K$ over the poset~$(P, \leq)$,
    \emph{the support of $s
    \in K$}, denoted $\supp_F(s)$, is defined as
    \[
        \supp_F(s):=\{u\in P| s\in F(u)\}.
    \]
    
    We call any $x \in P$ such that $x \in \supp_F(s)$ and $s \notin
    F(x')$ for any $x'<x$ in $P$ an \emph{entrance coordinate} for $s$. Note
    that $s$ can have multiple entrance coordinates.
\end{definition}

Some filtrations are ``nicer'' than others.  For example, one nice property is
if the support of
every simplex has a minimum (i.e., unique) \emph{entrance coordinate}.  If a
filtration
has this property, we call it \emph{one-critical}.  In this paper, we focus on slighly more general filtrations - those where each simplex has at most $k\in\N$ entrance coordinates.

\begin{definition}[one critical, $k$-critical, \cite{bauer2025fast}]
A bifiltration is called \emph{free or 1-critical} if the support
of every simplex is a principal upset, that is, the upward closure of a single element in the parameter space. More generally, a bifiltration is \emph{$k$-critical} if the support of a simplex is the
upward closure of at most $k$ elements.
\end{definition}

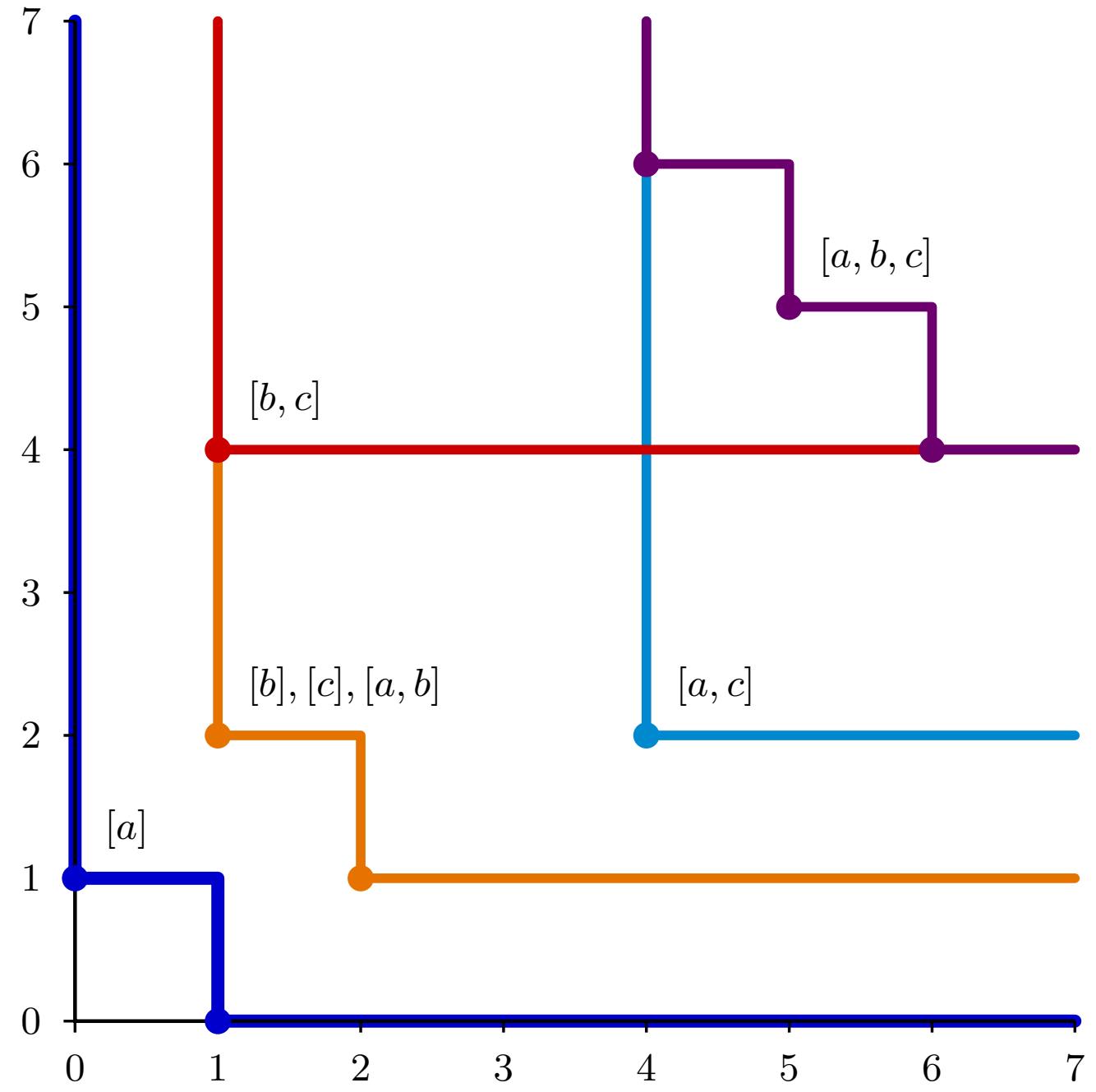
\begin{figure}[h!]
    \centering
     {\begin{tikzpicture}
        \input{figs/tixzfigures/biparameterspace_example}
    \end{tikzpicture}}
    \caption{This figure displays the entrance coordinates and entrance thresholds of the simplicies of the bifiltration from \Cref{fig:filtration}.  
    Note that this filtration is $3$-critical. 
    }\label{fig:entrancethresh}
\end{figure}

\begin{rems}
   When the parameter space $P$ of a $k$-critical filtration $F$ has $P\subset \R^d$, for a simplex $s\in\KK$, we refer to the boundary of $\supp_F(s)$ as the \emph{entrance threshold of $s$}, denoted $\partial\supp_F(s)$. 

    We sometimes refer to the \emph{support of a filtration $F$}, this is simply the union of the supports of all simplices in $\KK$.  
\end{rems}

\subsection{Persistence Modules}

Applying the homology functor to a filtration, we obtain a module of vector
spaces that we call a persistence module. More generally, a persistence module is defined as follows:

\begin{definition}[Persistence Module]
    Let $\K$ be a finite field and $(P,\leq)$ a poset.
    A \emph{persistence
    module} over~$(P,\leq)$ with coefficients in $\K$ is a
    functor $M \colon (P,\leq) \to \vec{\K}$, where $\vec{\K}$ denotes the
    category of $\K$-vector spaces.  In
    other words, $M$ assigns:
    \begin{itemize}
        \item For each point $\vecx \in P$, a graded $\K$-vector space $M_{\vecx}$.
        \item For each relation $\vecx \leq \vecy$ in $P$,
            a linear map $M_{\vecx, \vecy}: M_{\vecx} \rightarrow M_{\vecy}$
            such that $M_{\vecx,
            \vecx}$ is the identity map and $M_{\vecy, \vecz} \circ M_{\vecx,
            \vecy} = M_{\vecx, \vecz}$.
    \end{itemize}

    We often call $P$ the parameter space for $M$.  We denote the rank of the map
    $M_{\vecx,\vecy}$ as $\beta^{\vecx,\vecy}$.

    If $P=(\R^d, \preceq)$ or $P=[\veca,\vecb] \subset (\R^d, \preceq)$
    and $\dim(P)=d$,
    then we call~$M$ a~\emph{$d$-parameter persistence module}.  Note
    that $d$-parameter filtrations
    induce $d$-parameter persistence modules.  In the special case when $d=2$, we may also refer to $M$ as a bipersistence module.
\end{definition}


While persistence modules need not arise from a filtration, this is a common
setting.
To see how this works, consider a simplicial complex $K$
and a filtration $F: (P,\leq) \rightarrow 2^K$ of
nested subcomplexes of a simplicial complex~$K$.
Then, the persistence module \emph{corresponding to or induced by}~$F$ is the persistence
module~$M \colon (P,\leq) \to \vec{\K}$ obtained from applying the
homology functor.
This yields a persistence module
where~$M_{x} = H(F(x); \K)$
and, for $x \leq y$,  the map~$M_{x, y}$ is the linear map on the homology groups
induced from the inclusion~$F(x) \hookrightarrow F(y)$.
When~$\K$ is clear from context, we simply write~$H(F(x))$ for~$H(F(x); \K)$.
When $P=\Rbar:=\R\cup\infty$, we obtain the common construction of one-parameter persistence modules.

If a filtration $F$ is $k$-critical, then the induced persistence module is particularly nice.  One immediate consequence of $k$-criticality for a bi-filtration is that the associated bipersisence module is \emph{$S$-constructible.}  $S$-constructible modules are defined for $1$-parameter persistence modules in \cite{patel2018generalized}, and generalized in \cite{clause2025meta}.

\begin{definition}[$S$-constructible modules, \cite{clause2025meta}]
A persistence module $M:(\mathbb{R},\leq)\to \vec{\K}$ is \emph{constructible} if there exists a finite set $S:=
\{s_1 <\dots<s_n\}\subset\mathbb{R}$ such that:
\begin{itemize}
    \item For $a<s_1, M_a = 0$;
    \item For $s_i \leq a\leq b<s_{i+1}$, $M_{a,b}$ is an isomorphism;
    \item For $s_n \leq a\leq b$, $M_{a,b}$ is an isomorphism.
\end{itemize}

Similarly, a bipersistence module $M : (\mathbb{R}^2,\leq) \to \vec{\K}$ is constructible if there exists a finite set $S := \{s_1 <\dots<s_n\}\subset R$ such that:
\begin{itemize}
    \item If $x<s_1$ or $y<s_1$, then $M_{(x,y)} = 0$,
    \item For $s_i \leq x_1 \leq x_2 < s_{i+1}$ and $s_j \leq y_1 \leq y_2 < s_{j+1}, M_{(x_1,y_1),(x_2,y_2)}$ is an isomorphism,
    \item If $x_1 \geq s_n$ or $y_1 \geq s_n$ and $(x_1,y_1) \leq (x_2,y_2)$, then $M_{(x_1,y_1),(x_2,y_2)}$ is an
isomorphism.
\end{itemize}

In either case, such a module is $S$-constructible.
\end{definition}

In the case of a $k$-critical bifiltration $F$, the set $S$ may be derived as a subset of all $(x,y)$ coordinates for all entrance coordinates of simplices in $\KK$.  For example, the set $S$ for the fitration given in \Cref{fig:filtration} is $S=\{0,1,2,4,5,6\}.$

\begin{definition}[Persistence Diagram]
    A (one-parameter) persistence diagram~$D$ is a multiset of
    dimension-labeled points in the
    extended plane $\Rbar^2$.  More specifically, a persistence diagram~$D$ is a
    finite multiset of
    elements in~$\Rbar^2 \times \N$, where $\N=\{0,1,2, \ldots \}$.  In particular
    $(b,d,k) \in D$ represents a
    $k$-dimensional feature that is \emph{born} at
    time~$b$ and dies going into time $d$. \footnote{Here, we are using ``time'' to stand
    for the one-dimensional persistence parameter.}

    We sometimes denote as $D^k$ the subset of $D$ of points with dimension $k$ fixed.

\end{definition}

An equivalent way to visualize persistence diagrams is with
\emph{persistence barcodes}. A point $(b,d,k)$ in the persistence diagram
corresponds to one copy of the interval $[b,d)$ in the persistence barcode.
We use
barcodes to illustrate a few examples throughout this paper, and refer to the barcode of a one-parameter persistence module $M$ as $B(M)$.

\begin{definition}[Mapping Persistence Modules to Diagrams]\label{def:pmtodgm}
    Let~$M \colon \Rbar \to \vec{\K}$ be a graded persistence module, and $k \in \N$.

    Then, the
    \emph{$k$-dimension persistence diagram
    corresponding to or induced by~$M$}, denoted $D^k(M)$, is the
    multiset of points in $\Rbar^2$ where $(x,y) \in \Rbar^2$ has multiplicity:
    \begin{equation}\label{eq:multiplicity}
        \mu_k(x,y)
        = \lim_{\epsilon\to0^+} (\beta_k^{x,y-\epsilon}- \beta_k^{x,y})
        -(\beta_k^{x-\epsilon,y-\epsilon}-\beta_k^{x-\epsilon,y}).
    \end{equation}

    By construction, multiplicity for a point $(x,y)$ is $0$ if $y\leq x$.  We
    include $(x,y)$ in the persistence diagram if and only if it has positive
    multiplicity.

    The persistence diagram of $M$, denoted $D(M)$, is defined as the multi-set
    \[
        D(M)=\{(b,d,k)^{\mu_k(b,d)}|(b,d)\in D^k(M)\}.
    \]
\end{definition}

In \Cref{eq:multiplicity}, the first difference counts the number of homology
classes that are born at or before $x$ and die at $y$. The second difference
counts the number of homology classes that are born at or before $x-\epsilon$ and
d, i.e., at $y$. We subtract these two values so that we do not double count.

\begin{definition}[Space of Persistence Diagrams, Space of Persistence Barcodes]
    Let $\diagrams$ denote the set of all persistence diagrams (i.e.,
    finite multisets of points in $\Rbar^2 \times \N$).
    Let $d_{PD} \colon \diagrams \times \diagrams \to \R$
    be a metric on the set of persistence diagrams (e.g., the
    bottleneck distance between diagrams).
    The space of persistence diagrams is the metric space~$(\diagrams, d_{PD})$.
    The topology on this space is the metric topology.

    Equivalently, one can define the \emph{space of persistence barcodes}, $\mathcal{B}$, also using the bottleneck distance as its metric.
\end{definition}

The generalization of a barcode for $k-d$ persistence modules with $k>1$ is not a trivial task; as such, we look to find other possible invariants of $k-d$ persistence modules.  One such invariant for bipersistence modules is the \emph{fibered barcode}, a summary obtained by fibering parameter space with lines of positive slope, in order to obtain a collection of $1-d$ persistence modules.

Originally defined in \cite{cerri2013betti} (and then named in \cite{lesnick2015interactive}), we show here the definition of fibered barcode from \cite{lesnick2025fast}.

\begin{definition}[Fibered Barcode \cite{lesnick2025fast}]
Let $M \colon (P,\leq) \to \vec{\K}$ be a bipersistence module. Let $L\subset\R^2$ be an affine line with positive slope.  Let $M_L : (L,\leq)\to \vec{\K}$ denote the restriction
of $M$ to $L$.  Note that, as the restriction of the partial order on $\R^2$ to $L$ is a total order, $M_L$ is a one-parameter persistence module with a well-defined barcode $B(M_L)$. 

This allows one to define, for a given module $M$, a map from $\mathcal{L}$, the set of affine lines in $\R^2$ with positive, finite slope, to the space of persistence barcodes $\mathcal{B}$. This map is called \emph{the fibered barcode of $M$}, and is denoted as $F(M)$, and given by $F(M)(L)=\mathcal{B}(M_L)$.

\end{definition}

\section{The Directed Parameter Space Complex}\label{sec:dpsc}

\subsection{The Directed Parameter Space Complex}






In what follows, we will fix the following set up for our filtrations and modules.  We fix the indexing poset $P$ as $P=\R^2_+=\{(x,y)\in\R^2|x\geq0\textrm{ and }y\geq0\}$. We also assume that $F:(\R^2_+,\leq)\to2^{\KK}$ is a $k$-critical bifiltration of a simplicial complex $\KK$ with$\{c_1,\dots,c_m\}\subset \R^2_+$ the set of entrance points in $\R^2_+$ for all simplices in $\KK$. Finally, we assume that all entrance points $c_i$ have integer coordinates.

Since $F$ is $k$-critical, one can then define $S=\{s_0<s_1<\dots<s_l\}\subset\mathbb{R}$ as the totally ordered set of all possible $(x,y)$ coordinates of entrance points of simplices in $F$.  Note that it is  the case that all entrance coordinates for simplices for $F$ are contained in the cube $C=[(0,0),(s_l,s_l)].$  

It must also be the case that $F((x_1,x_2))=\KK$ for all $x_1,x_2>s_l$.  Moreover, \begin{itemize}
    \item $F(x_1,x_2)=F(x_1,s_l)$ for all $x_2>s_l$ and
    \item $F(x_1,x_2)=F(s_l,x_2)$ for all $x_1>s_l$. 
\end{itemize} 

Additionally, there is a natural Euclidean cubical structure on $C$, given by the elementary cubes in $C$ with integer coordinates. If $\mathcal{C}$ denotes the set of open cubes defined by points in $C\cap\Z^2$, then $(C,\mathcal{C},\di{C})$ is a directed Euclidean cubical complex. We use this structure to define the \emph{directed parameter space complex}.

\begin{definition}[Directed Parameter Space Complex]\label{def:dpsc}
Let $F:(\R^2_+,\leq)\to2^\KK$ be a $k$-critical bifiltration of a simplicial complex $\KK$, with associated persistence module $M:(\R^2_+,\leq)\to\vec{\K}$ obtained by applying the homology functor.

Using the steps listed above, we define $(C,\mathcal{C},\di{C})$ to be the associated \emph{directed parameter space complex of $F$ (or $M$).} 
\end{definition}

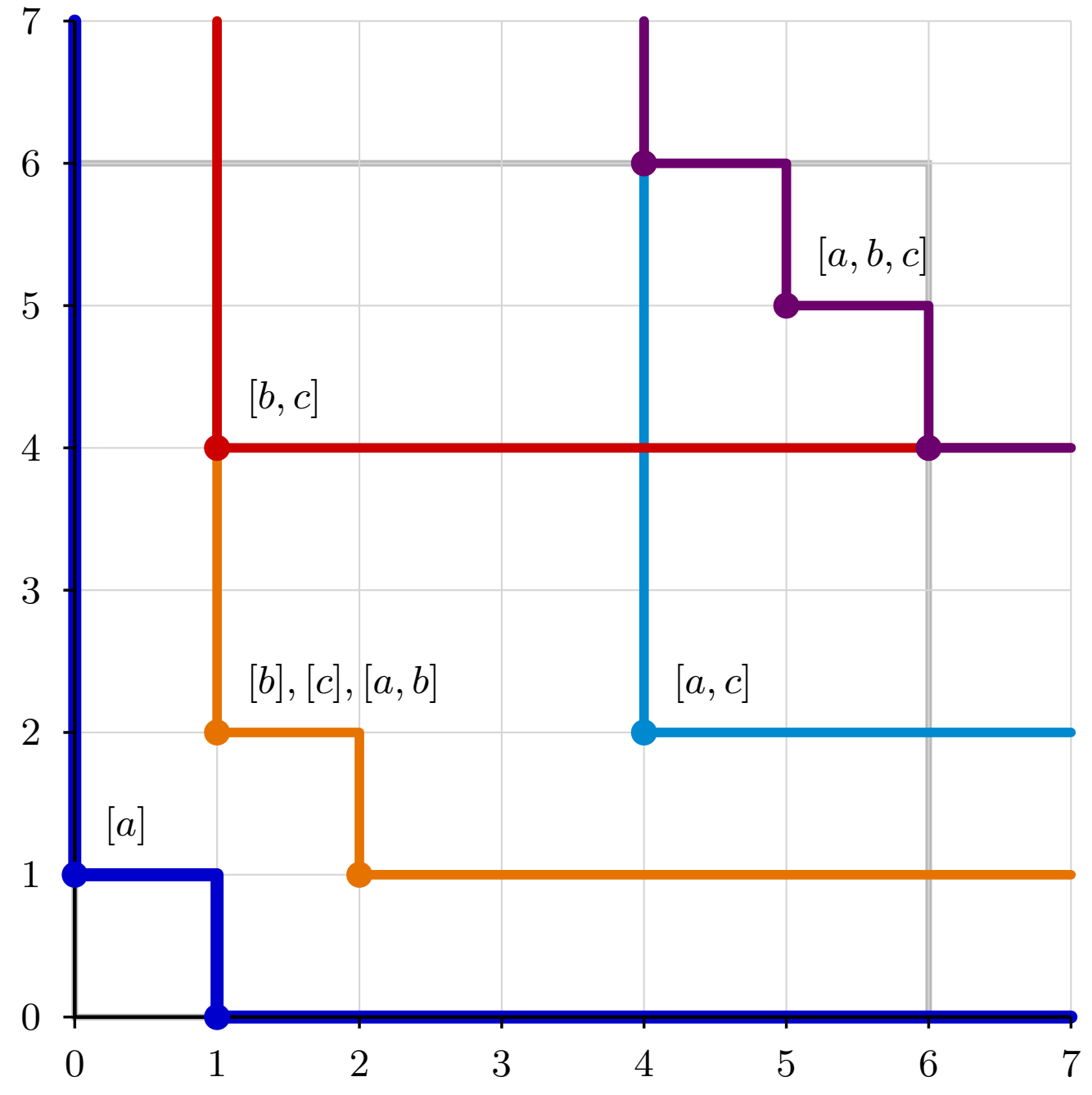
\begin{figure}[h!]
    \centering
     {\begin{tikzpicture}
        \input{figs/tixzfigures/directedPS_example}
    \end{tikzpicture}}
    \caption{The directed parameter space complex for the filtration in \Cref{fig:filtration}.  The boundary of $C=[(0,0),(6,6)]$ is shown in light grey.}\label{fig:dpsc_example}
\end{figure}

For example, for the filtration given in \Cref{fig:filtration}, $C=[(0,0),(6,6)].$  We illustrate the cubical structure for $C$ in \Cref{fig:dpsc_example}

\subsection{Modules along directed paths}

Recall, for a filtration $F$ defined over a poset $P$, we may obtain restriced filtrations $F_{P'}$ for subposets $P'\subset P$.  In particular, if $F$ is a $k$-critical filtration with indexing poset $P=\R^2_+$ with associated directed parameter space complex $(C,\mathcal{C},\di{C})$, each path $\alpha\in\di{C}$ induces a one-parameter filtration $F_\alpha$, a $1-d$ persistence module $M_\alpha$ with associated persistence diagram denoted $D(\alpha)$, and persistence barcode denoted $B(\alpha)$.  Note that $M_\alpha$ is a \emph{sub-module} of $M$, the bipersistence module associated to $F$.\footnote{Such sub-modules are studied by restricting to
lines with positive slope in $\R^2$ ~\cite{lesnick2015interactive} and the corresponding software, RIVET~\cite{rivet}.}

Observe that for each $\alpha\in\di{C}$, $F_\alpha$ is determined completely by the order that $\alpha$ intersects the entrance thresholds of simplices in $\KK$ in  the following way.  First, the following definition is given to provide language and notation to track how entrance coordinates for simplicies in $F$ relate to entrance coordinates for simplicies in $F_\alpha.$

\begin{definition}[Upset, Push]

For a point $c\in \R^2$, define the {\em upset of $c$ in $\R^2$},
denoted $U_c$, as
$$U_c:=\{x\in\R^2 |c\preceq x\}.$$

Given a non-decreasing path $\alpha:[0,1]\to \R^2$ such that $\textrm{im}(\alpha)\cap U_c\neq\varnothing$, we define the \emph{push of $c$ along $\alpha$} as 
\[\push_\alpha(c)=\min_{x\in \textrm{im}(\alpha)\cap U_c}x.\]

Note that $\push_\alpha(c)$ is well defined, unique, and contained in $\textrm{im}(\alpha)\cap U_c$.  Finally, we denote $p_\alpha(c)$ to be the value in $[0,1]$ for which
$\push_\alpha(c)=\alpha(p_\alpha(c))$. 
\end{definition} 

When $F$ is $k$-critical with indexing poset $P=\R^2_+$ with associated directed parameter space complex $(C,\mathcal{C},\di{C})$, for each simplex $\sigma\in\KK$, there are a finite
number of incomparable entrance coordinates $\{e^1,\dots,e^k\}\subset\R^2_+$ for
$\sigma$, and
\[
    \supp_F(\sigma)=\bigcup_{j=1}^kU_{e^j}.
\]

This means there is at least one entrance coordinate $e^*\in\{e^1,\dots,e^l\}$ for
which $p_\alpha(e^*)$ is minimal in the set
$\{p_\alpha(e^1),\dots,p_\alpha(e^l)\}$; $e^\ast$ is unique unless $\alpha$ passes through the least upper bound of two entrance coordinates $e_i$ and $e_j$, where the least upper bound of two points $u,v\in\RR^2$ is defined as $\lub(u,v)=(\max(u_1,v_1),\max(u_2,v_2))$. 

We may think of
$\push_\alpha(e^*)$ as the point in $\R^2_+$ where $\sigma$
enters the filtration $F_\alpha$ along $\im(\alpha)$, and may sometimes denote $\push_\alpha(e^*)=\push_\alpha(\sigma)$ and $p_\alpha(e^*)=p_\alpha(\sigma)$, to reinforce the fact that it is the entrance of $\sigma$ along $\alpha$.  

\begin{figure}[h!]
    \centering
     {\begin{tikzpicture}
        \input{figs/tixzfigures/path_in_dpsc_ex}
        \end{tikzpicture}\subcaption[]{A dipath $\alpha$ in the directed parameter space complex for the filtration in \Cref{fig:filtration}. Simplex labels for the entrance thresholds are omitted. The push of entrance values for each simplex are indicated along $\alpha$, but only one is labeled to avoid figure overcrowding.  Note that a minimal set of entrance coordinates determining the filtration on $\alpha$ is the set $\{(0,1),(1,2), (1,4),(4,2),(5,5)\}$.}\label{fig:path_example}
        \begin{tikzpicture}
        \input{figs/tixzfigures/path_dpsc_filtration}
    \end{tikzpicture}\subcaption[]{The $1$ parameter filtration for $\alpha$.}}
    \caption{A path and its restricted filtration.}
\end{figure}
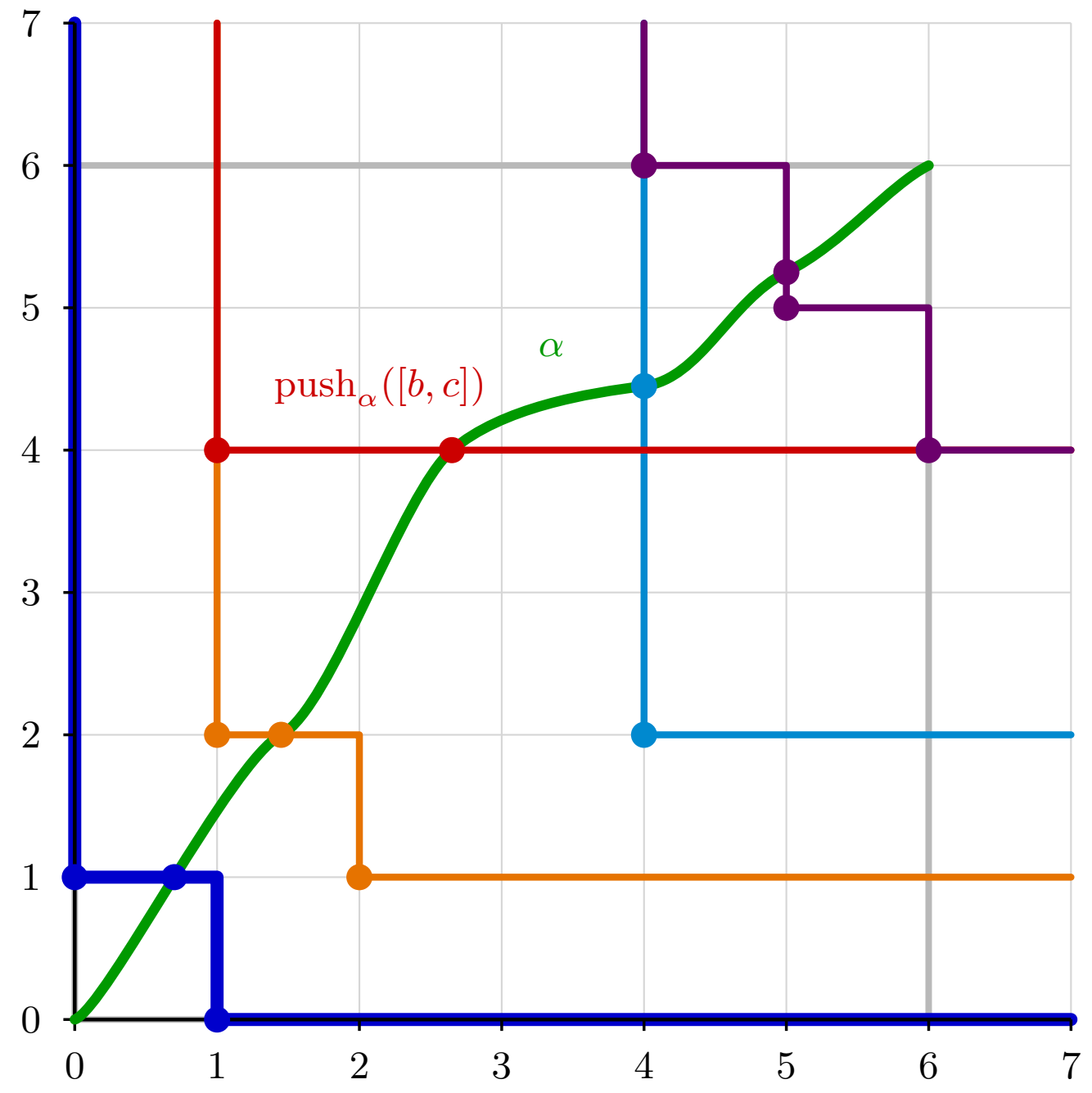
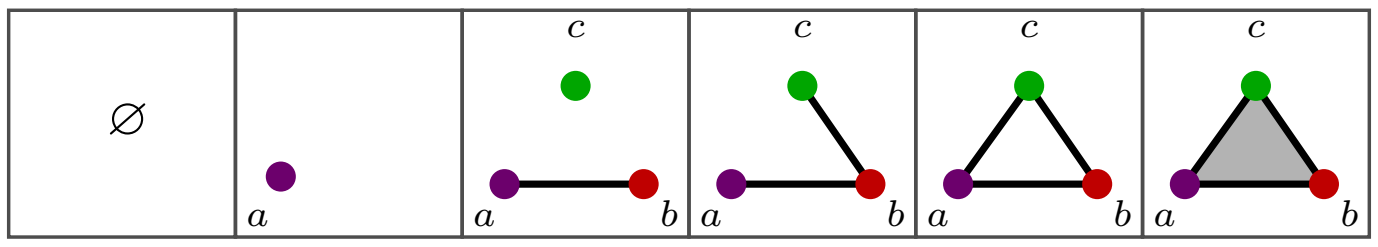

It is important to note that, for a dipath $\alpha,$ simplices can only enter $F_\alpha$ at a parameter value determined (via the push) by an entrance coordinate of that simplex in $F$.  As changes in homology in the corresponding module $M_\alpha$ can only occur with the entrance of simplex, the criticial homological parameter values of $M_\alpha$ are can be derived from the set of entrance coordinates of simplices along $\alpha$.  

In other words, given a path $\alpha$, let $0\leq t_1<t_2<\dots<t_j\leq1$ be the set of entrance coordinates of simplices in $F_\alpha$.  Then each $t_i=p_\alpha(\sigma_i)$ for some simplex $\sigma_i$, and as points in $\R^2_+$, $\alpha(t_i)=\push_\alpha(\sigma_i)$.  Furthermore, for all $i$, any values $t_i\leq s<s'<t_{i+1}$ have $F_s=F_{s'}$ and therefore $M_s=M_{s'}$.  

We summarize this observation in a corollary.

\begin{corollary}\label{cor:entrance_homology}
    Let $F:(\R^2_+,\leq)\to2^\KK$ be a $k$-critical bifiltration of a simplicial complex $\KK$ with associated persistence module $M:(\R^2_+,\leq)\to\vec{\K}$ obtained by applying the homology functor and 
    associated directed parameter space complex $(C,\mathcal{C},\di{C})$.
    
    For each let $\alpha\in\di{C}$, $M_\alpha$ is $S-$constructible with $S\subset\{p_\alpha(\sigma_i)|\sigma_i\in\KK\}.$
\end{corollary}

\section{Equivalence of Persistence Modules}\label{sec:moduleequivalence}
 In this section, we introduce two types of 1-d persistence modules.  First, we define combinatorial equivalence; two 1-d persistence modules are combinatorially equivalent when there is a order-preserving bijection between the birth-death coordinates of the associated persistence diagram.  While useful, combinatorial equivalence can be hard to determine along paths in the parameter space complex.  This motivates our second equivalence - simplicial equivalence.  Two filtrations of a simplicial complex $\KK$ are said to be simplicially equivalent if simplices enter in the same relative order.  In the following sections, we provide definitions, examples, and some properties for both types of equivalence.

\subsection{Combinatorial Equivalence of Persistence Diagrams}

In this section, we define combinatorial equivalence between persistence
diagrams. 

A persistence diagram uses a finite set of birth and death coordinates.  Let
$\eta$ denote the number of distinct values (including $\infty$). If
we forget these values and just consider the poset of coordinates,
we have a bijection between coordinates in the persistence diagram and the
(unique) totally ordered set of~$\eta$ elements.
More specifically:

\begin{definition}[Index Function and Combinatorial Maps]
    Let $D$ be a persistence diagram and $\coords(D) \subseteq \Rbar$ be the set
    of all birth and death coordinates in $D$. The \emph{index function}
    $\idx_D: \coords(D) \to \Nbar:=\N\cup\infty$ maps a birth or death coordinate in $D$ to
    its index in the total order of the
    elements of $\coords(D)$.

    Let $\rho_D:D\to \Nbar^2 \times \N$ be the map that
    takes a point~$(a,b,k)\in D$
    to~$\rho_D(a,b,k):=(\idx_D(a),\idx_D(b),k)$.
\end{definition}

\begin{example} Consider the diagram

\[
    D = \{ (0, \infty, 0), (1,5,1), (6,8,1) \}.
\]

Then, the total order of the elements of $\coords(D)$ is:
$0<1<5<6<7<\infty$, which is in one-to-one correspondence with the totally ordered poset:
$1<2<3<4<5<\infty$.
So,~$\idx_D(0)=1$, $\idx_D(1)=2$, and so forth.

\label{ex:coords}
\end{example} 

The idea is that we want two persistence diagrams to be combinatorially equivalent when the relative ordering between birth-death pairs is the same.
Notice that, at a minimum, for two persistence diagrams to be
combinatorially equivalent they must have the same number of birth-death pairs
for each possible dimension $k$.

\begin{definition}[Combinatorial and Combinatorially Equivalent Persistence
    Diagrams]
    Let $A$ be a one-dimensional persistence diagram.

    We define the \emph{combinatorial map} $\rho:\diagrams \to \diagrams$ by
    \[
    \rho(A):\{(\idx_A(b),\idx_A(d),k)|(b,d,k)\in A\}.
    \]

    The \emph{combinatorial (persistence) diagram of $A$} is defined as $\rho(A)$ 
    and encodes the relative orderings of all births and deaths in
    $A$. The \emph{combinatorial (persistence barcode) of $A$} is defined as the barcode of the combinatorial diagram $\rho(A)$. Equivalently, it is the image of the barcode $B(A)$ under mapping the birth/death coordinates using the index function $\idx_A$.

    Two diagrams $A$ and $B$ are \emph{combinatorially equivalent} if and only if
    $\rho(A)=\rho(B)$; in fact, they are equivalent if and only if
    $[\rho(A)]^k=[\rho(B)]^k$ for all $k$.  Moreover, within each equivalence
    class of combinatorially equivalent persistence diagrams, there exists a
    unique diagram $D^*=\rho(D^*)$.
    \label{def:ce}
\end{definition}

    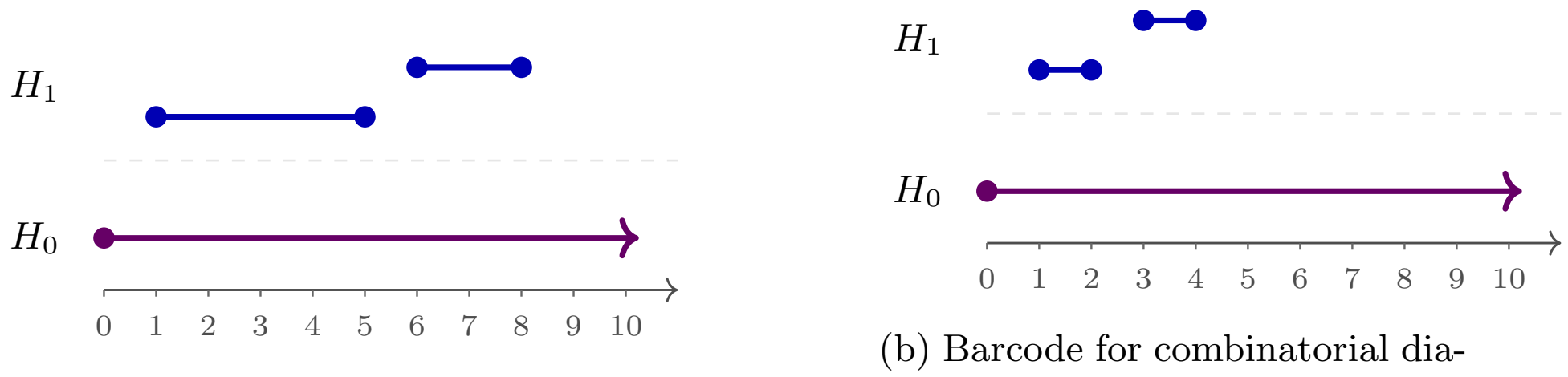
\begin{figure}[htp]
        \centering
        \begin{subfigure}[b]{0.4\textwidth}
            \centering
            \begin{tikzpicture}
                \input{figs/tixzfigures/barcode_A}
            \end{tikzpicture}
            \caption{Barcode for diagram $A$}
            \label{subfig:barcodes-original}
        \end{subfigure}
        \hfill
        \begin{subfigure}[b]{0.4\textwidth}
            \centering
            \begin{tikzpicture}
                \input{figs/tixzfigures/barcode_rhoA}
            \end{tikzpicture}
            \caption{Barcode for combinatorial diagram of $A$, $\rho(A)$}
            \label{subfig:barcodes-combinatorial}
        \end{subfigure}
        \caption{Barcode $\to$ Combinatorial Barcode. In
            \subfigref{barcodes}{original}, we visualize the barcode given in~\exref{coords} where $A=\{(0,\infty, 0), (1,5,1), (6,8,1)\}$. In \subfigref{barcodes}{combinatorial}, we visualize the combinatorial barcode $\rho(A)=\{(0,\infty, 0), (1,2,1), (3,4,1)\}$ which captures the ordering of the birth-death pairs but forgets the lengths of the bars.}
        \label{fig:barcodes}
    \end{figure}

Next we show that indeed, combinatorial equivalence is an equivalence relation on the
set of persistence diagrams.

\begin{lemma}[CE is an Equivalence Relation]\label{lem:equiv}
    Combinatorial equivalence is an equivalence
    relation on the set of one-dimensional persistence diagrams.
\end{lemma}
\begin{proof}
By definition,   two diagrams $A$ and $B$ are combinatorially equivalent if and only if
    $\rho(A)=\rho(B)$. Because equality is reflexive, symmetric and transitive, so is combinatorial equivalence.
\end{proof}

\begin{figure}[htp]
        \centering
        \begin{subfigure}[b]{0.4\textwidth}
            \centering
            \begin{tikzpicture}
                \input{figs/tixzfigures/barcode_A}
            \end{tikzpicture}
            \caption{Barcode for diagram $A$}
            \label{subfig:nonCEbarcode-1}
        \end{subfigure}
        \hfill
        \begin{subfigure}[b]{0.4\textwidth}
            \centering
            \begin{tikzpicture}
                \input{figs/tixzfigures/barcode_D}
            \end{tikzpicture}
            \caption{Barcode for diagram $D$}
            \label{subfig:nonCEbarcode-2}
        \end{subfigure}
        \caption{Non Combinatorial Equivalent Barcodes. We visualize the barcodes given in~\exref{CEdiagrams} that are not combinatorially equivalent. In
            \subfigref{nonCEbarcode}{1}, $A=\{(0,\infty, 0), (1,5,1), (6,8,1)\}$. In \subfigref{nonCEbarcode}{2}, $D=\{(0,\infty, 0), (1,5,1), (5,8,1)\}$. Since the (strict) ordering of the birth-death events are not preserved, they are not combinatorially equivalent. Hence, the combinatorial barcodes are not equal: $\rho(A)\neq\rho(D)$.}
        \label{fig:nonCEbarcode}
    \end{figure}
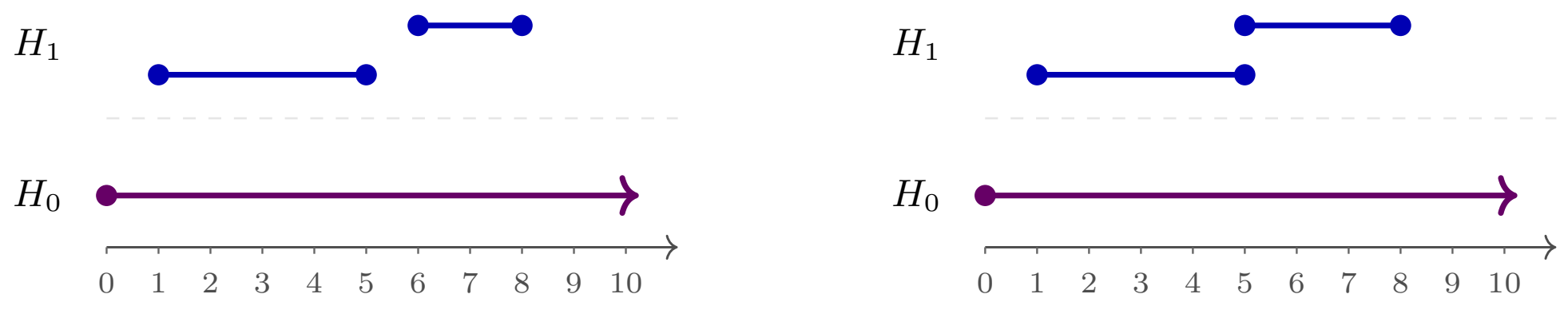

We
provide an example and non-example  of combinatorially equivalent diagrams
next:

\begin{example}[Combinatorial Equivalent Diagrams] \label{ex:CEdiagrams}
    Consider the following persistence diagrams
    \begin{equation}
        A=\left\{(0,\infty,0),(1,5,1),(6,8,1)\right\}
        \text{ and }
        B=\left\{(0.5,\infty,0),(1,5,1),(7,8,1)\right\}.
    \end{equation}
    Each has one
    zero-dimensional point and two one-dimensional points.
    The total order of~$\coords(A)$ is~$0<1<5<6<8<\infty$ and
    the total order of~$\coords(B)$ is~$0.5<1<5<7<8<\infty$.
    These two
    diagrams are combinatorially equivalent, and the equivalence class representative
    is $\rho(A)=\rho(B)=\left\{(0,\infty,0),(1,2,1),(3,4,1)\right\} $.

    However, if we consider another persistence diagram $$D=\{(0,\infty,0),(1,5,1),(5,8,1)\},$$
    then~$A$ and $D$ are not combinatorial equivalent. The total order on the birth and
    death times in $D$ is~$0<1<5<8<\infty$.
    At a minimum, \defref{ce}  requires
    the diagrams to have the same number of distinct coordinates, which is not the case
    here.
    The equivalence class representative  of $D$
    is $\rho(C)=\left\{(0,\infty,0),(1,2,1),(2,3,1)\right\} $. See~\figref{CEbarcode} for visualizing combinatorially equivalent barcodes and~\figref{nonCEbarcode} for visualizing barcodes that are not combinatorially equivalent.
\end{example}

    \begin{figure}[htp]
        \centering
        \begin{subfigure}[b]{0.4\textwidth}
            \centering
            \begin{tikzpicture}
                \input{figs/tixzfigures/barcode_A}
            \end{tikzpicture}
            \caption{Barcode for diagram $A$}
            \label{subfig:CEbarcode-1}
        \end{subfigure}
        \hfill
        \begin{subfigure}[b]{0.4\textwidth}
            \centering
            \begin{tikzpicture}
                \input{figs/tixzfigures/barcode_B}
            \end{tikzpicture}
            \caption{Barcode for diagram $B$}
            \label{subfig:CEbarcode-2}
        \end{subfigure}
        \caption{Combinatorial Equivalent Barcodes. We visualize the barcodes given in~\exref{CEdiagrams}. In
            \subfigref{CEbarcode}{1}, $A=\{(0,\infty, 0), (1,5,1), (6,8,1)\}$. In \subfigref{CEbarcode}{2}, $B=\{(0,\infty, 0), (1,5,1), (7,8,1)\}$. Since the ordering of the birth-death events are preserved, they are combinatorially equivalent. Hence, the combinatorial barcodes are equal: $\rho(A)=\rho(B)$.}
        \label{fig:CEbarcode}
    \end{figure}
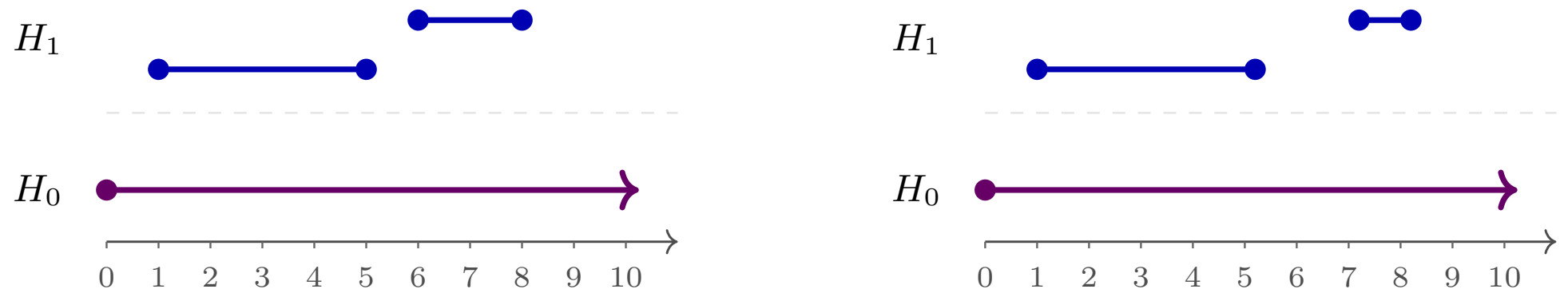

Using~\lemref{equiv}, we define the following quotient space.

\begin{definition}[Space of Combinatorial Diagrams]
    The space of persistence diagrams modulo combinatorial equivalence with the
    quotient topology is the \emph{space of combinatorial diagrams}, and is
    denoted $\cdiagrams$.
\end{definition}

\begin{definition}[Space of Combinatorial Barcodes]
    The space of persistence barcodes modulo combinatorial equivalence with the
    quotient topology is the \emph{space of combinatorial barcodes}, and is
    denoted $\mathcal{CB}$.
\end{definition}

We can now represent a point in the space of combinatorial diagrams by the unique representative given by the map $\rho$. In the following remarks we distinguish between combinatorial filtrations and dipaths, as well as paths and dipaths.

\begin{rems}[Combinatorially Equivalent Filtrations and Dipaths]
    We say that two filtrations $F$ and $F'$ are {\em combinatorially
        equivalent} if they induce combinatorially equivalent persistence
    diagrams. 
    
    We also say that two dipaths $\alpha$ and $\beta$ are
    combinatorially equivalent if $D(\alpha)$ and $D(\beta)$ are combinatorially
    equivalent.
\end{rems}


\begin{definition}[Combinatorial Fibered Barcode]
Let $M \colon (P,\leq) \to \vec{\K}$ be a bipersistence module and let $F(M)$ be the fibered barcode of $M$. We can extend the map $F(M)$ to the space of combinatorial barcodes $\mathcal{CB}$ by simply composing with the map $\rho$ from \Cref{def:ce}.  We call this map the \emph{combinatorial fibered barcode of $M$}, denoted $CF(M)$.

\end{definition}

\subsection{Simplical Equivalence of Filtrations}\label{sec:SE}
Tracking homological changes locally (with only a small perturbation of a path $\alpha$) is not straight forward. As such, we define an additional relationship between filtrations, \emph{simplicial equivalence}.  Although this is a finer equivalence than combinatorial equivalence, it is much more practical to partition the set of directed paths for a filtration via simplical equivalence than combinatorial equivalence.  Additionally, simplicial equivalence implies combinatorial equivalence.

\begin{definition}[Simplicial Equivalence]
Let $(P,\leq)$ be a totally ordered poset, and let $F$ be a filtration of a simplicial complex $\KK$ defined over $(P,\leq)$.  Let $\{e_1,\dots,e_k\}\subset P$ be the set of entrance coordinates for simplices in $\KK$ in $F$.

We define the \emph{filtration index function}, denoted $\textrm{idx}_F$, as the function which maps a simplex entrance coordinate to its index in the total order of the the elements of $\{e_1,\dots,e_k\}$.

For each filtration $F$, we can then define a filtration $\tilde{F}:\N\to 2^\KK$ as the filtration where every simplex $\sigma$ with entrance coordinate $e_i$ in $F$ has entrance coordinate $\textrm{idx}_F(e_i)$ in $\tilde{F}$.

We say that two filtrations $F,G$ are \emph{simplicially equivalent} if $\tilde{F}=\tilde{G}.$
\end{definition}

By definition, simplical equivalence is an equivalence relation. 
\begin{corollary}[SE is an Equivalence Relation]\label{lem:seequiv}
    Simplicial equivalence is an equivalence
    relation on the set of one-dimensional filtrations on a fixed simplicial complex $\KK$.
    
\end{corollary}

We obtain the following corollary by observing that if two filtrations  are simplicially equivalent, then they must have combinatorially equivalent diagrams.

\begin{corollary}\label{cor:reparam_CE}
Let $F:(\R^2_+,\leq)\to2^\KK$ be a $k$-critical bifiltration with associated directed parameter space complex $(C,\mathcal{C},\di{C})$, and let $\alpha,\beta\in\di{C}$.  If $F_\alpha,F_\beta$ are simplicially equivalent, then $\rho(D(F_\alpha))=\rho(D(F_\beta))$.
\end{corollary}

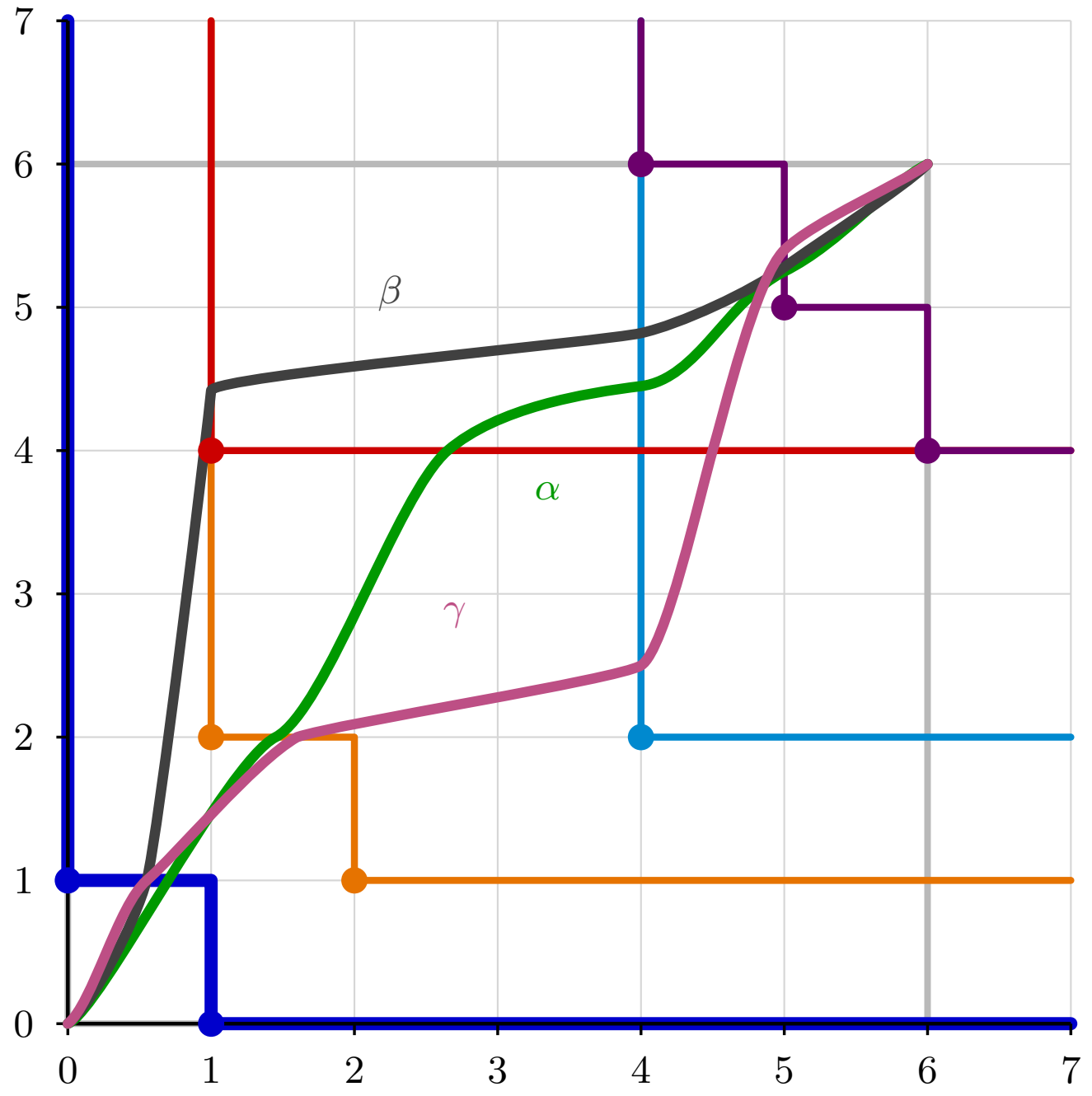
\begin{figure}[h!]
    \centering
    \begin{tikzpicture}
        \input{figs/tixzfigures/incomparablelemma_ex}
    \end{tikzpicture}
    \caption{Examples of paths which are not simplicially equivalent. When comparing $\alpha$ and $\beta$, we see \Cref{lem:incomparable}(\ref{item:incomparable}) is satisfied, and when comparing $\alpha$ and $\gamma$ we see that \Cref{lem:incomparable}(\ref{item:incomparable}) is satisfied.}
    \label{fig:incomparable}
\end{figure}

The converse of \Cref{cor:reparam_CE} is not true; it is possible for the diagrams along two paths to be combinatorially equivalent but not simplicially equivalent.  For example, the paths $\alpha$ and $\gamma$ in \Cref{fig:incomparable} are combinatorially but not simplically equivalent, due to the edges $[a,c]$ and $[b,c]$ entering in different orders.

The following lemma gives sufficient conditions for two paths to not be simplicialy equivalent; if two filtrations along paths $\alpha$ and $\beta$ (with the same starting and ending points) are not simplically equivalent, then there must be nontrivial intersection of the boundaries of the supports of two simplicies $\sigma$ and $\tau$ in the region bounded by $\alpha$ and $\beta.$

\begin{lemma} \label{lem:threshold}

Let $F:(\R^2_+,\leq)\to2^\KK$ be a $k$-critical bifiltration with associated directed parameter space complex $(C,\mathcal{C},\di{C})$, let $\alpha,\beta\in\di{C}$ with $\alpha(0)=\beta(0)$ and $\alpha(1)=\beta(1)$, and let $R$ be the subset of $\RR^2$ bounded by $\alpha$ and $\beta$.  If the filtrations $F_\alpha,F_\beta$ are not simplicially equivalent, then there exist two distinct simplicies $\sigma$ and $\tau$ for which $\sigma,\tau\notin F_{\alpha(0)}$, $\sigma,\tau\in F_{\alpha(1)}$, and $R\cap\partial[\supp_F(\sigma)]\cap\partial[\supp_F(\tau)]\neq\varnothing$.

\end{lemma}
\begin{proof}


Assume $F_\alpha$ and $F_\beta$ are not simplicially equivalent.  If it was the case that for all pairs of simplices $\sigma$ and $\tau$ for which $\sigma,\tau\notin F_{\alpha(0)}$, $\sigma,\tau\in F_{\alpha(1)}$, it is the case that $R\cap\partial[\supp_F(\sigma)]\cap\partial[\supp_F(\tau)]=\varnothing$, then  either 
\[R\cap\supp_F(\sigma)\subset R\cap\textrm{int}(\supp_F(\tau)) \qquad\textrm{or}\qquad R\cap\supp_F(\tau)\subset R\cap\textrm{int}(\supp_F(\sigma)) \]  This would imply that $p_\alpha(\sigma)<p_\sigma(\tau)$ and  $p_\beta(\sigma)<p_\beta(\tau)$ or $p_\alpha(\sigma)>p_\sigma(\tau)$ and  $p_\beta(\sigma)>p_\beta(\tau)$, i.e., $\sigma$ appears before $\tau$ along both $\alpha$ and $\beta$, or $\tau$ appears before $\sigma$ along both $\alpha$ and $\beta$.  Since this is true for all pairs, this implies $F_\alpha$ and $F_\beta$ are simplically equivalent, a contradiction.
\end{proof}

In addition to \Cref{lem:threshold}, \Cref{lem:incomparable} provides necessary and suffient conditions on simplicial equivalence along paths with the same start and end points. This is useful for the following sections, as it provides a tangible description of what occurs when two paths are not simplicially equivalent.

\enlargethispage{3\baselineskip}
\begin{lemma} \label{lem:incomparable}
Let $F:(\R^2_+,\leq)\to2^\KK$ be a $k$-critical bifiltration with associated directed parameter space complex $(C,\mathcal{C},\di{C})$, and let $\alpha,\beta\in\di{C}$ with $\alpha(0)=\beta(0)$ and $\alpha(1)=\beta(1)$.  The filtrations $F_\alpha,F_\beta$ are not simplicially equivalent if and only if there exist two entrance coordinates $c_\sigma\neq c_\tau$, for $\sigma,\tau\in\KK$ respectively, such that, without loss of generality exchanging $\alpha$ and $\beta$ or $\sigma$ and $\tau$, either
\begin{enumerate}
    \item \label{item:incomparable}$c_\sigma$ and $ c_\tau$ are incomparable, $\push_\alpha(c_\sigma)<\push_\alpha(c_\tau)$ and $\push_\beta(c_\sigma)>\push_\beta(c_\tau),$ or
    \item  \label{item:comparable}$c_\sigma$ and $ c_\tau$ are equal in exactly one coordinate and $\push_\alpha(c_\sigma)<\push_\alpha(c_\tau)$ and $\push_\beta(c_\sigma)=\push_\beta(c_\tau)$
\end{enumerate}
\end{lemma}

\begin{proof}
Note, two filtrations $F_\alpha$ and $F_\beta$ are simplicially equivalent if and only if, for all pairs of simplices $\sigma,\tau\in\KK$, either 
\begin{itemize}
    \item $p_\alpha(\sigma)<p_\alpha(\tau)$ and $p_\beta(\sigma)<p_\beta(\tau)$, or 
    \item $p_\alpha(\sigma)=p_\alpha(\tau)$ and $p_\beta(\sigma)=p_\beta(\tau)$, or 
    \item $p_\alpha(\sigma)>p_\alpha(\tau)$ and $p_\beta(\sigma)>p_\beta(\tau)$.
\end{itemize}

Thus, $F_\alpha$ and $F_\beta$ are not simplicially equivalent if and only if there exists a pair of simplices $\sigma,\tau\in\KK$ such that 

\[p_\alpha(\sigma)<p_\alpha(\tau)\qquad\textrm{  and }\qquad p_\beta(\sigma)\geq p_\beta(\tau).\]

These equalities hold if and only if there are entrance coordinates $c_\sigma$ and $c_\tau$, 

\[\push_\alpha(c_\sigma)<\push_\alpha(c_\tau)\qquad\textrm{  and }\qquad \push_\beta(c_\sigma)\geq\push_\beta(c_\tau)\] 
as desired.

Note that $c_\sigma^1<c_\tau^1$ and $c_\sigma^2<c_\tau^2$ if and only $\push_\gamma(c_\sigma)<\push_\gamma(c_\tau)$ for any monotonically increasing path $\gamma$. Therefore, any two entrance coordinates $c_\sigma$ and $c_\tau$ which satisfy the above inequalities must be incomparable or equal in exactly one coordinate.

So, the result holds.
\end{proof}

Finally, we show that, if $\alpha$ and $\beta$ are not simplicially equivalent, then at least one pair of simplicies $\sigma$ and $\tau$ obtained from \Cref{lem:threshold} satisfies the conditions of the simplicies obtained from \Cref{lem:incomparable}.

\begin{lemma} \label{lem:samesimplices}

Let $F:(\R^2_+,\leq)\to2^\KK$ be a $k$-critical bifiltration with associated directed parameter space complex $(C,\mathcal{C},\di{C})$, let $\alpha,\beta\in\di{C}$ with $\alpha(0)=\beta(0)$ and $\alpha(1)=\beta(1)$, and let $R$ be the subset of $\RR^2$ bounded by $\alpha$ and $\beta$.  If the filtrations $F_\alpha,F_\beta$ are not simplicially equivalent, then at least one pair of simplicies $\sigma$ and $\tau$ obtained in \Cref{lem:threshold} have entrance coordinates $c_\sigma\neq c_\tau$ which satisfy the conditions of \Cref{lem:incomparable}.
\end{lemma}

\begin{proof}
    Suppose $F_\alpha,F_\beta$ are not simplicially equivalent, and $\sigma$ and $\tau$ are simplicies such that $\sigma,\tau\notin F_{\alpha(0)}$, $\sigma,\tau\in F_{\alpha(1)}$, and $R\cap\partial[\supp_F(\sigma)]\cap\partial[\supp_F(\tau)]\neq\varnothing$, as obtained by \Cref{lem:threshold}.  

Note that the set $R\cap\partial[\supp_F(\sigma)]\cap\partial[\supp_F(\tau)]$ is the disjoint union of at least one of the following types of sets: 

 \begin{itemize}
     \item points $\vecp$, where $\vecp=\textrm{l.u.b}(c_\sigma,c_\tau)$ where $c_\sigma$ and $c_\tau$ are incomparable entrance coordinates of $\sigma$ and $\tau$ respectively, or
     \item segments defined as the edge between an entrance coordinate $c_\sigma$ and either $\push_\alpha(\sigma)$ or $\push_\beta(\sigma)$, where  $\textrm{l.u.b}(c_\sigma,c_\tau)=c_\sigma$, or
      \item segments defined as the edge between an entrance coordinate $c_\tau$ and either $\push_\alpha(\tau)$ or $\push_\beta(\tau)$, where  $\textrm{l.u.b}(c_\sigma,c_\tau)=c_\tau$.
 \end{itemize}

If there is at least one point $\vecp$ in the intersection, then we have two incomparable entrance coordinates $c_\sigma$ and $c_\tau$.  Otherwise, the intersection must have at least one segment defined by entrance coordinates $c_\sigma$ and $c_\tau$ which are equal in exactly one coordinate.

If it was the case that for all of the pairs obtained by \Cref{lem:threshold}, none of pairs of entrance coordinates generating elements of $R\cap\partial[\supp_F(\sigma)]\cap\partial[\supp_F(\tau)]$ satisfy the conditions of either (1) or (2) in \Cref{lem:incomparable}, then we would have simplical equivalence of $\alpha$ and $\beta$, which is a contradiction.

\end{proof}
 \section{The Extended Fibered Barcode}\label{sec:fb}

Given a fixed bi-filtration $F:(\R^2_+,\leq)\to2^\KK$ with associated directed parameter space complex $(C,\mathcal{C},\di{C})$, which possible filtrations or modules may be obtained along paths in $\di{C}$?  To address this question, we make the following observations and definitions about the modules obtained this way.

First, we note that there may be dipaths $\alpha$ in $\di{C}$ for which $F_\alpha$ does not start with the empty set, nor end with the final complex $\KK$. For example, a straight line path from $(1.5,4)$ to $(4.5,4)$ in our example filtration \Cref{fig:filtration} does not include the simplex $[a,b,c]$.

As opposed to considering all such restricted modules, we instead choose to only consider paths that bisect $C$ as a line with positive slope would - namely, paths which intersect both the past boundary $\partial^-(C)$ and the future boundary $\partial^+(C)$. Denote this set of paths as 
\[\dipaths{C}{l}{u}:=\displaystyle\bigcup_{\vecx\in \partial^-(C), \\ \vecy\in \partial^+(C)}\dipaths{C}{\vecx}{\vecy};\]

\noindent i.e, all paths that start in the past boundary of $C$ and terminate in the future boundary. 

\begin{figure}[h]
    \centering
     {\begin{tikzpicture}
        \input{figs/tixzfigures/efb_example}
        \end{tikzpicture}\caption{Lines with positive slope and the corresponding simplicially equivalent paths in $\dipaths{C}{\vecx}{(l,l)}$, from \Cref{prop:fiberedtoextended}.  The path $\alpha_{L_i}$ is simplicially equivalent to $L_i$ for $i=1,2,3$}
        }\label{fig:casesefb}
\end{figure}
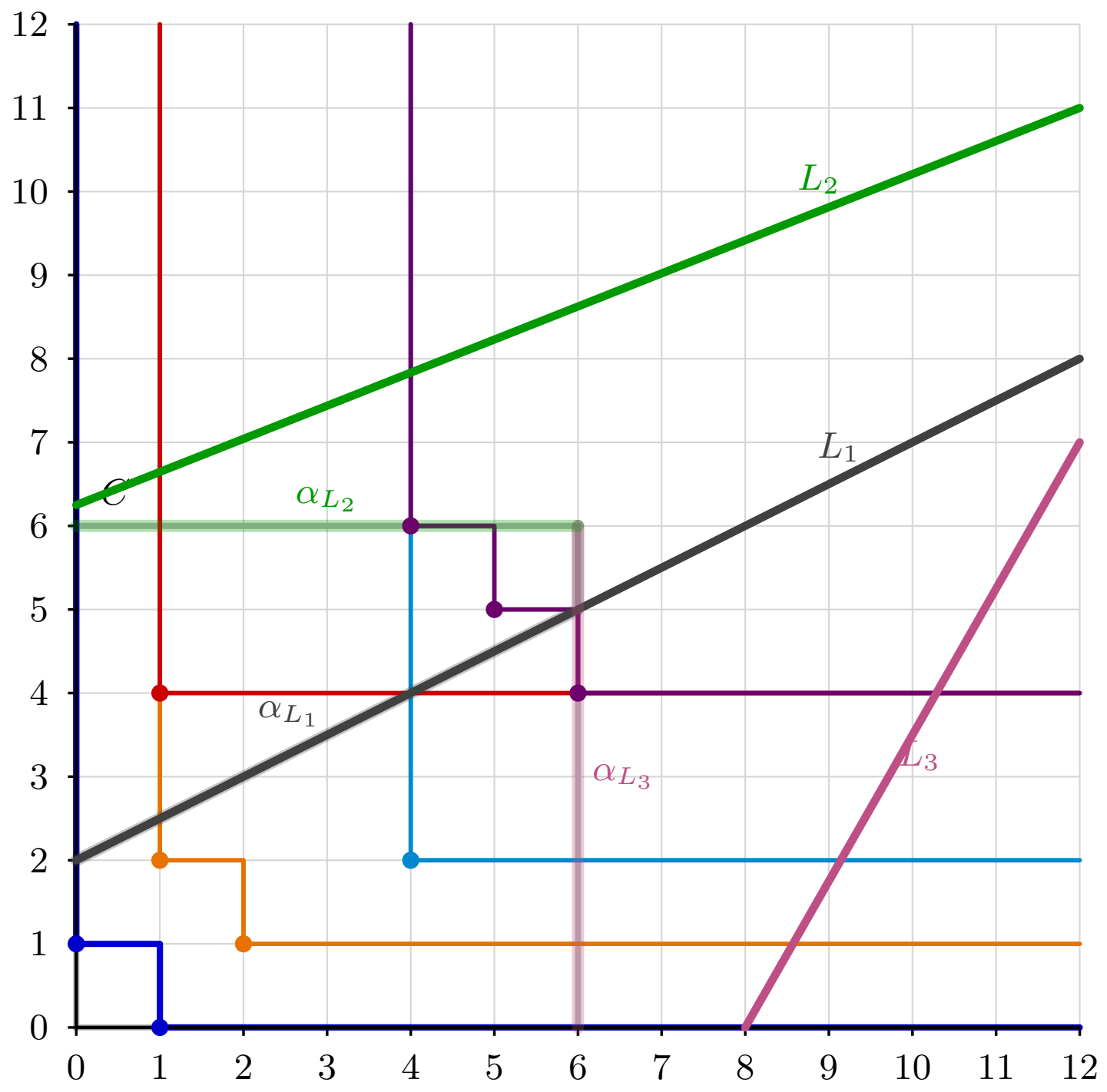

In \Cref{prop:fiberedtoextended}, we describe a correspondence between filtrations along paths in the directed parameter space complex to filtrations obtained along lines with positive slope. In particular, for any line $L$ with positive slope, there is always a path $\alpha$ in the directed parameter space complex with the property that $F_\alpha$ and $F_L$ are simplicially equivalent, as in \Cref{fig:casesefb}.

\enlargethispage{2\baselineskip}
 
\begin{proposition}\label{prop:fiberedtoextended}
     Let $F:(\R^2_+,\leq)\to2^\KK$ be a $k$-critical bi-filtration with associated bipersistence module $M$ and directed parameter space complex $(C, \mathcal{C}, \di{C})$.  Then for all lines $L$ with positive slope, there exists $\vecx\in \partial^-(C)$ and $\alpha\in \dipaths{C}{\vecx}{(l,l)}$ where $F_\alpha$ is simplicially equivalent to $F_L$.
\end{proposition}

\begin{proof}

If $L\cap C=\varnothing$, $F_L$ is simplicially equivalent to either:

\begin{itemize}
    \item the filtration along the straight horizontal line path from $(0,l)$ to $(l,l)$, if there exists a point $(0,r)\in L$ such that $r>l$, or
     \item the filtration along the straight vertical line path from $(l,0)$ to $(l,l)$, if there exists a point $(r,0)\in L$ such that $r>l$,
\end{itemize}
both of which are in $\dipaths{C}{\vecx}{(l,l)}$.

If $L\cap C\neq\varnothing$, then there exists $\vecx\in \partial^-(C)$ such that $\vecx\in L$.  If this is the case, then the filtration $F_L$ is a simplicially equivalent to $F_\alpha$ where $\alpha(0)=\vecx$, $\alpha$ follows along $L$ until the point $L\cap \partial^+(C)$ and then follows the future boundary monotonically until $(l,l)$. 

\end{proof}

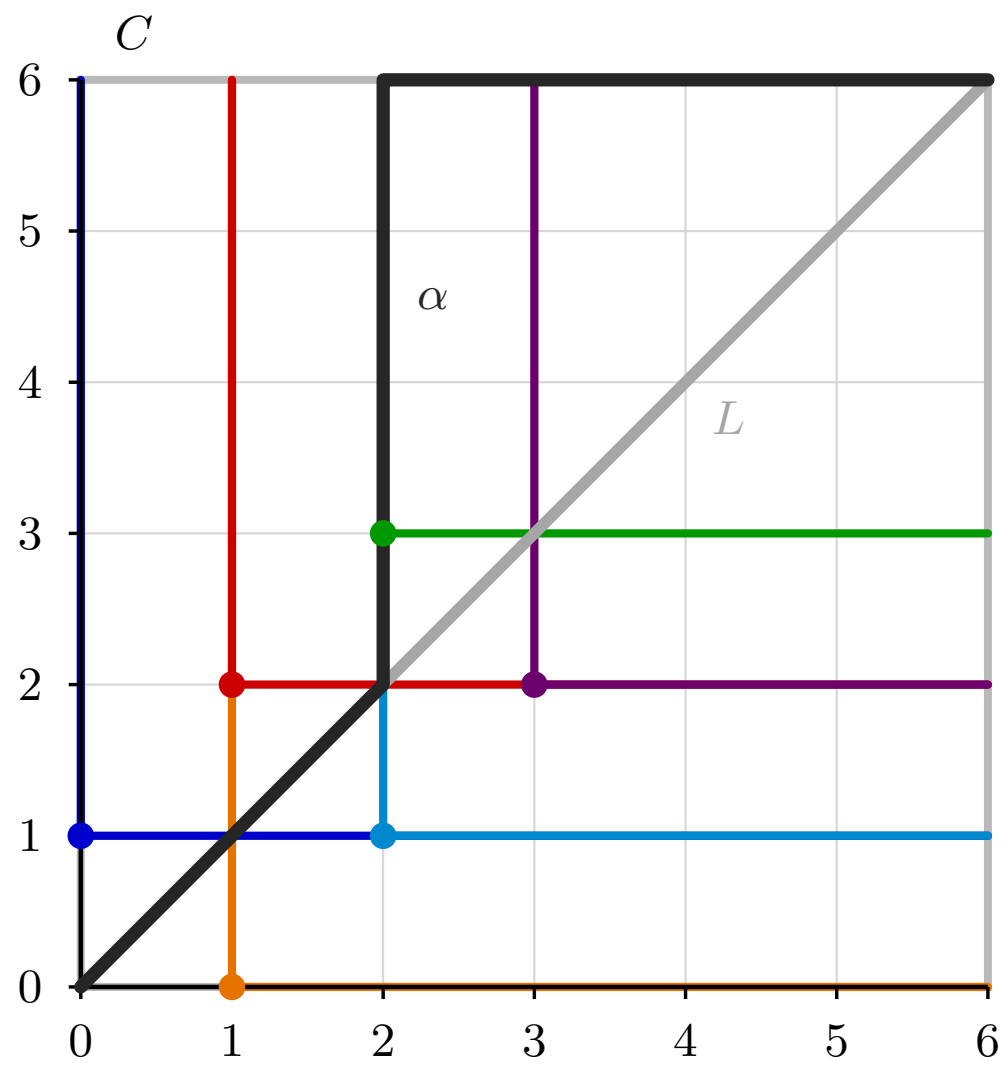
\begin{figure}[h!]
    \centering
     {\begin{tikzpicture}
        \input{figs/tixzfigures/six_thresholds_alpha}
    \end{tikzpicture}}
    \caption{An example of a path $\alpha$ for which there is no line with positive slope which has a filtration which is simplicially equivalent to $\alpha$. If an entrance threshold indicates the entrance of a different simplex, then any line where which passes through the least upper bound of the entrance coordinates for orange and blue, and the least upper bound of cyan and red, must also pass through the least upper bound of green and purple, indicating these simplicies enter at the same point along $L$.  However, there exist paths, such as $\alpha$, for which green may enter before purple.}\label{fig:pathnotin_fb}
\end{figure}

\enlargethispage{2\baselineskip}
It is then clear that by looking at the set of filtrations obtained along paths in $\dipaths{C}{l}{u}$, we obtain a set of filtrations that, modulo simplical equivalence, at least contains all possible filtrations obtainable by fibering with lines with positive slope. In many cases, there will be filtrations along paths in $\dipaths{C}{l}{u}$ that are not obtainable on a line with positive slope (see \Cref{fig:pathnotin_fb}.). This motivates our next definition, the extended fibered barcode, to include the possible new modules and barcodes obtained along paths in $\dipaths{C}{l}{u}$.

 \begin{definition}[Extended Fibered Barcode, Extended Combinatorial Fibered Barcode.]\label{def:efb}
     Let $C=[(0,0),(l,l)]\subset\R^2$, $M \colon (C,\leq) \to \vec{\K}$ be a bipersistence module, and $\alpha\in\dipaths{C}{l}{u}$.  Recall, $M_\alpha : (\alpha,\leq)\to \vec{\K}$ denotes the restriction of $M$ to $\alpha$. 

This allows one to define, for a given module $M$, a map from $\dipaths{C}{l}{u}$, to the space of persistence barcodes $\mathcal{B}$. This map is called \emph{the extended fibered barcode of $M$}, and is denoted as $EF(M)$, and given by $EF(M)(\alpha)=\mathcal{B}(M_\alpha)$.
When composed with the combinatorial map $\rho$, we obtain the \emph{extended combinatorial fibered barcode of $M$}, denoted $ECF(M):=\rho\circ EF(M)$.
 \end{definition}

Simplicially equivalent filtrations will induce combinatorially equivalent diagrams. But we are able to say even more than that; given $L$ and $\alpha$ as obtained in \Cref{prop:fiberedtoextended}, there will always exist a non decreasing-reparameterization $r$ such that $D(F_L)=D(F_{\alpha\circ r})$. 
So, it is the case that the image of original fibered barcode is contained in the image of the extended fibered barcode.

\begin{corollary}\label{cor:fb_in_efb}
  Let $F$ be a $k$-critical bi-filtration with associated bipersistence module $M$, let $F(M)$ be the fibered barcode of $M$ and let $EF(M)$ be the extended fibered barcode of $M$. Then $\textrm{im}(F(M))\subset\textrm{im}(EF(M))$.
    
\end{corollary}

We also have the following statements about the combinatorial versions of the fibered and extended fibered barcode, \Cref{lem:ecf_finite} and \Cref{lem:fibered_finite}.  Essentially, \Cref{lem:ecf_finite} says that since we are dealing with a finite simplicial complex, there's only a finite number of ways that you can order the simplices. In other words, the size of $im(ECF)$ is less than $2^m$ where $m$ is the number of simplices in $\KK$.

\begin{lemma}\label{lem:ecf_finite}

Let $F$ be a $k$-critical bi-filtration, with associated bipersistence module $M$ and let $ECF(M)$ be the extended combinatorial fibered barcode of $M$.

Then $\textrm{im}(ECF(M))\subset\mathcal{CB}$ is a finite set.

\end{lemma}

\begin{proof}
Let $F:(\R^2_+,\leq)\to2^{\KK}$ be a $k$-critical bifiltration of a simplicial complex $\KK$.

For each $\vecx\in \partial^-(C)$, a path $\alpha\in \dipaths{C}{\vecx}{(l,l)}$ induces $F_\alpha$, a filtration of $\KK$.  As there are finitely many simplicies in $\KK$, there are finitely many simplicial equivalence classes of filtrations of $\KK$, determined by the total ordering of entrance coordinates of simplices along $\alpha$.

So, the set of filtrations obtained along paths $\alpha\in \dipaths{C}{\vecx}{(l,l)}$ has a finite number of equivalence classes under simplicial equivalence.  Each simplicial equivalence class maps to exactly one combinatorial equivalence class in the space of barcodes.  Thus, the image of $ECF(M)$ in $\mathcal{CB}$ is a finite set.
\end{proof}

As \Cref{cor:fb_in_efb} states that the image of the fibered barcode is contained in the image of the extended fibered barcode, we obtain the following corollary.

\begin{corollary}\label{lem:fibered_finite}

Let $F$ be a $k$-critical bi-filtration, with associated bipersistence module $M$ and let $CF(M)$ be the combinatorial fibered barcode of $M$.

Then $\textrm{im}(CF(M))\subset\mathcal{CB}$ is a finite set.
\end{corollary}

\section{The Entrance Complex}\label{sec:ec}

For a given $k$-critical filtration $F:(\R^2_+,\leq)\to2^{\KK}$ and associated persistence module $M$, we have an associated directed Euclidean cubical complex $(C,\mathcal{C},\di{C})$, and we will use the directed structure to analyze properties of the extended fibered barcode of $M$.  However, as $C$ is a  cube in $\R^2$, all paths in $\di{C}$ and therefore $\dipaths{C}{l}{u}$ are dihomotopy equivalent. 

We aim to identify a cubical subcomplex of $C$ for which the following is true: dihomotopy equivalence of paths 
implies simplicial equivalence of barcodes along those paths. As we have the information of the entrance coordinates of simplices for $F$, and homological changes along a path $\alpha$ can only occur at points $\textrm{push}_\alpha(c_i)$ for $c_i\in E$ (\Cref{cor:entrance_homology}), we will construct a subcomplex which omits certain $2$-cubes. These will be the ones whose future boundary intersects more than one entrance threshold for more than one simplex.  This, as we will prove later, prevents the existence of a dihomotopy that connects two dipaths for which those simplices enter in opposing (total) orders.

\begin{definition}[Entrance Subcomplex $(E,\mathcal{E},\di{E})$]\label{def:ec}

  Let  $F:(\R^2_+,\leq)\to2^{\KK}$ be a $k$-critical bifiltration with associated bipersistence module $M$, and ~$(C,\mathcal{C},\di{C})$ be the associated directed Euclidean cubical complex.

  We call a cube $c\in\mathcal{C}$ a \emph{mixed entrance $2$-cube} if there exists at least two simplices $\sigma_1,\sigma_2\in\KK$ such that $\partial\textrm{supp}_F(\sigma_i)\cap \partial^+(\overline{|c|})\neq\varnothing$ for $i=1,2$ and $\partial\textrm{supp}_F(\sigma_1)\cap \partial^+(\overline{|c|})\neq\partial\textrm{supp}_F(\sigma_1)\cap \partial^+(\overline{|c|})$.

The \emph{entrance subcomplex for $F$}, denoted $(E,\mathcal{E},\di{E})$ is the directed Euclidean subcomplex of $(C,\mathcal{C},\di{C})$ defined by:
\begin{itemize}
    \item $\mathcal{E}_0=\mathcal{C}_0$, $\mathcal{E}_1=\mathcal{C}_1$ and $\mathcal{E}_2$ contains only the $2$-cubes $\sigma\in \mathcal{C}$ that are not mixed-entrance $2$-cubes,
    \item $E=|\mathcal{E}|$, and
    \item $\di{E}$ is the subset of $\di{C}$ with paths contained in $E$.
\end{itemize}
\end{definition}

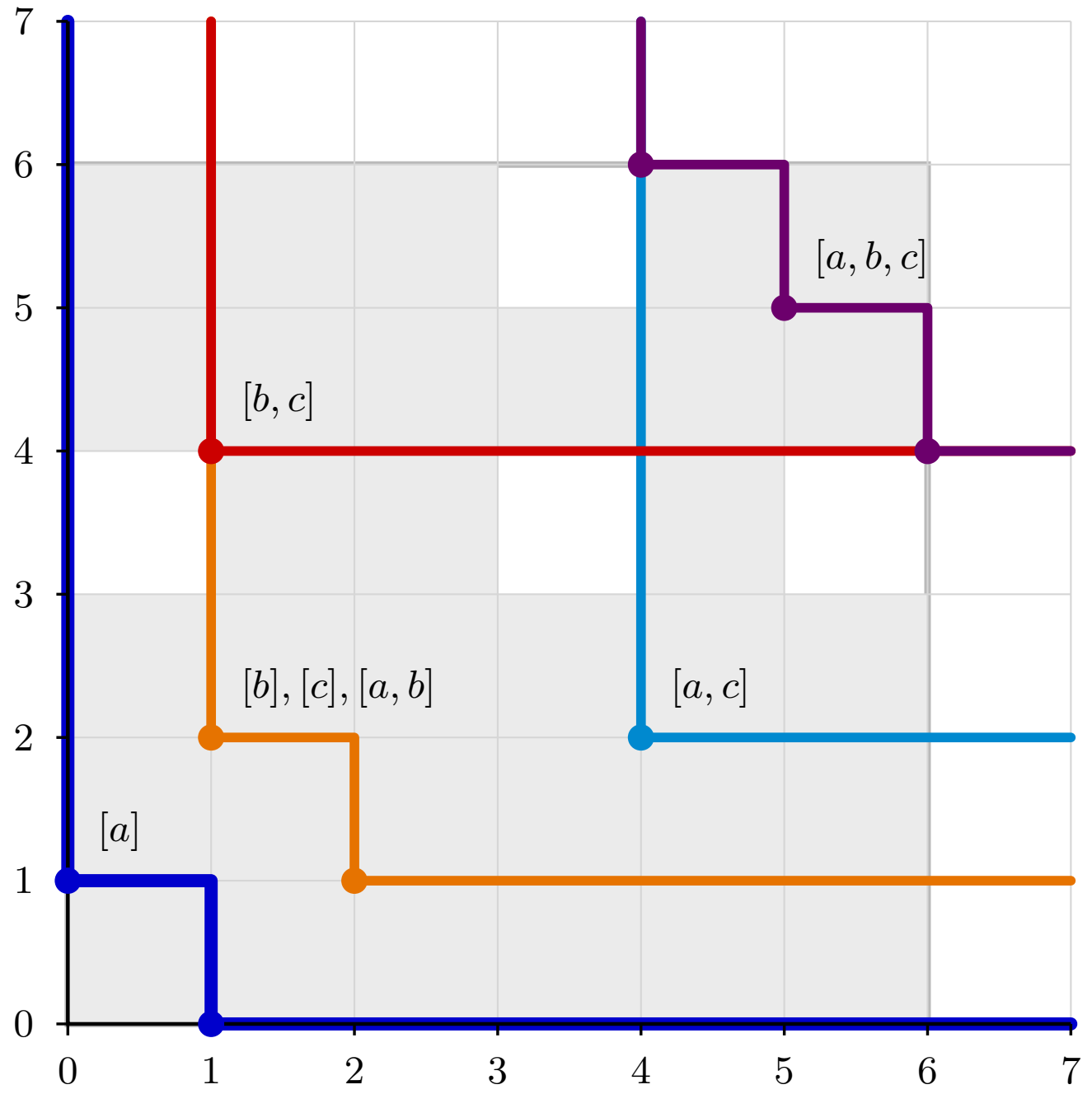
\begin{figure}[h]
    \centering
     {\begin{tikzpicture}
        \input{figs/tixzfigures/entrance_subcomplex}
        \end{tikzpicture}\caption{Entrance subcomplex for the directed parameter space complex of \Cref{fig:dpsc_example}, shown in grey. Mixed cubes are $[(0,3),(1,4)],[(3,3),(4,4)],[(5,3),(6,4)],$ and $[(3,5),(4,6)]$.}
        }\label{fig:entrancesubcomplex}
\end{figure}

In summary, the mixed cubes are areas where simplices could enter in different orders and the entrance complex is  $C=[ \veczero,(l,l)] $ with (the interior of) mixed cubes removed. Intuitively, a mixed cube $c\in\mathcal{C}$ is a cube for which $\dipaths{c}{l}{u}$ contains paths where at least two simplices $\sigma_1$ and $\sigma_2$ enter in different relative orders.  We remove them from the entrance complex $E$ to create dihomotopy obstructions between paths that are not simplicially equivalent.

 \begin{samepage}
In \Cref{lem:kcrit_dihomotopy}, we show that for a $k$-critical bifiltration, dihomotopy equivalence in the Entrance Complex implies simplicial equivalence. 

 \begin{theorem}\label{lem:kcrit_dihomotopy}
Let $F$ be a $k$-critical bifiltration with entrance complex $(E,\mathcal{E},\di{E})$, let $\vecx\prec\vecy\in E$ and let $\alpha,\beta\in \dipaths{E}{\vecx}{\vecy}$ with $\alpha(0)=\beta(0)$ and $\alpha(1)=\beta(1)$.  If there is a dihomotopy $H$ between $\alpha$ and $\beta$ with $H_t\in\dipaths{E}{\vecx}{\vecy}$ for all $t\in[0,1]$, then $\alpha$ and $\beta$ are simplicially equivalent. 
 \end{theorem}
 \begin{proof}
     We will proceed by contradiction and suppose $\alpha$ and $\beta$ are not simplicially equivalent. 
 By \Cref{lem:samesimplices}, there exist two simplicies $\sigma$ and $\tau$ for which $\sigma,\tau\notin F_{\alpha(0)}$, $\sigma,\tau\in F_{\alpha(1)}$, and $R\cap\partial[\supp_F(\sigma)]\cap\partial[\supp_F(\tau)]\neq\varnothing$, where $R$ is the region in $\RR^2$ bounded by $\alpha$ and $\beta,$ and two entrance coordinates $c_\sigma$, $c_\tau$ which satisfy the conditions of \Cref{lem:incomparable}.

     Let $\vecp$ be the least upper bound of $c_\sigma$ and $c_\tau$; note that $\vecp\in R$. We next show that cube given by $C^*:=[\vecp-\bf{1},\vecp]$ is mixed entrance. In the case that $c_\sigma$ and $c_\tau$ are incomparable as in \Cref{lem:incomparable} \ref{item:incomparable}, the intersection of $C^*$ with the entrance threshold of one simplex is $[(p_1-1,p_2),\vecp]$ and with the threshold of the other simplex is $[(p_1,p_2-1),\vecp]$, which makes $C^\ast$ mixed entrance. 
    
    Now consider the case that $c_\sigma$ and $c_\tau$ are comparable as in \Cref{lem:incomparable} \ref{item:comparable}. Since $\alpha$ and $\beta$ are not simplicially equivalent, without loss of generality, $\push_\alpha(\sigma)<\push_\alpha(\tau)$ and $\push_\beta(\sigma)=\push_\beta(\tau)$. The intersection of $C^*$ with the entrance threshold of one simplex is $[(p_1-1,p_2),\vecp]$ or $[(p_1,p_2-1),\vecp]$ and with the threshold of the other simplex is $\vecp$. This also implies that  $C^*$ is  a mixed entrance cube. 


    Now, consider the region $R$ bounded by the curves $\alpha$ and $\beta$.  
    Note, this region can be partitioned into pieces determined by which path $\alpha$ or $\beta$ is above or below, or if the paths are equal. We can write $R=R_1\cup\dots\cup R_k$ where for each $R_i$, either $\alpha$ is above, below or equal to $\beta$ in the region $R_i$, and $R_i$ contains exactly two points $\vecz_{i-1},\vecz_i\in\alpha\cap\beta$ with $\vecz_{i-1}\leq r\leq\vecz_i$ for all $r\in R_i$. 

    Since $\vecp\in R$, it must be contained in exactly one region $R^\ast\in\{R_1,\dots,R_k\}$.  Moreover, it cannot be the case that $R^\ast$ is a region where $\alpha$ and $\beta$ coincide, as this would imply equality for the pushes of both simplices on both paths.  So, without loss of generality, assume that $R^\ast$ is one where $\alpha$ is above and $\beta$ is below.  

Now we show that int$(C^\ast)\cap R^\ast\neq\varnothing$. 
If $\vecp$ lies neither on $\alpha$ nor $\beta$, then it is in the interior of $R^\ast$, and thus int$(C^\ast)\cap R^\ast\neq\varnothing$. It cannot be the case that $\vecp$ lies on both $\alpha$ and $\beta$.

If $\vecp$ is on $\alpha$, then $\alpha$ intersects $C^\ast$ and at least one of the cubes $[(p_1-1,p_2),\vecp]$, $[(p_1,p_2-1),\vecp]$, or int$(C^\ast)$.  If $\alpha$ intersects int$(C^\ast)$ or $[(p_1-1,p_2),\vecp]$ non trivially, as $\vecp$ does not lie on $\beta$ and $\beta$ is below $\alpha$, it must be the case that int$(C^\ast)\cap R^\ast\neq\varnothing$.  It cannot be the case that $\alpha\cap C^\ast\subseteq[(p_1,p_2-1),\vecp]$, as this would imply that $\tau$ enters before $\sigma$ along $\alpha$. A similar argument shows that int$(C^\ast)\cap R^\ast\neq\varnothing$ when $\vecp$ is strictly to the right of $\alpha$ and on $\beta$.

 Finally, because $E\subset \R^2 $, it is true that any dihomotopy between $\alpha$ and $\beta$ must contain the regions $R_1,\dots,R_k$ bounded by $\alpha$ and $\beta$ in its image  \footnote{To see that this is true, note that any dihomotopy between $\alpha$ and $\beta$ is, if direction is forgotten, a path homotopy with fixed endpoints in $\R^2$.  If such a homotopy does not cover the regions between the paths, then following \cite{munkres2000topology}  one could puncture $\R^2$ and have a path homotopy between two non homotopic paths in punctured $\R^2$.} Thus, $R^\ast\subset\textrm{im}(H)$ and int$(C^\ast)\cap R^\ast\neq\varnothing$ together imply that int$(C^\ast)\cap \textrm{im}(H)\neq\varnothing$, giving the contradiction required. 
 
 \end{proof}

 \end{samepage}

 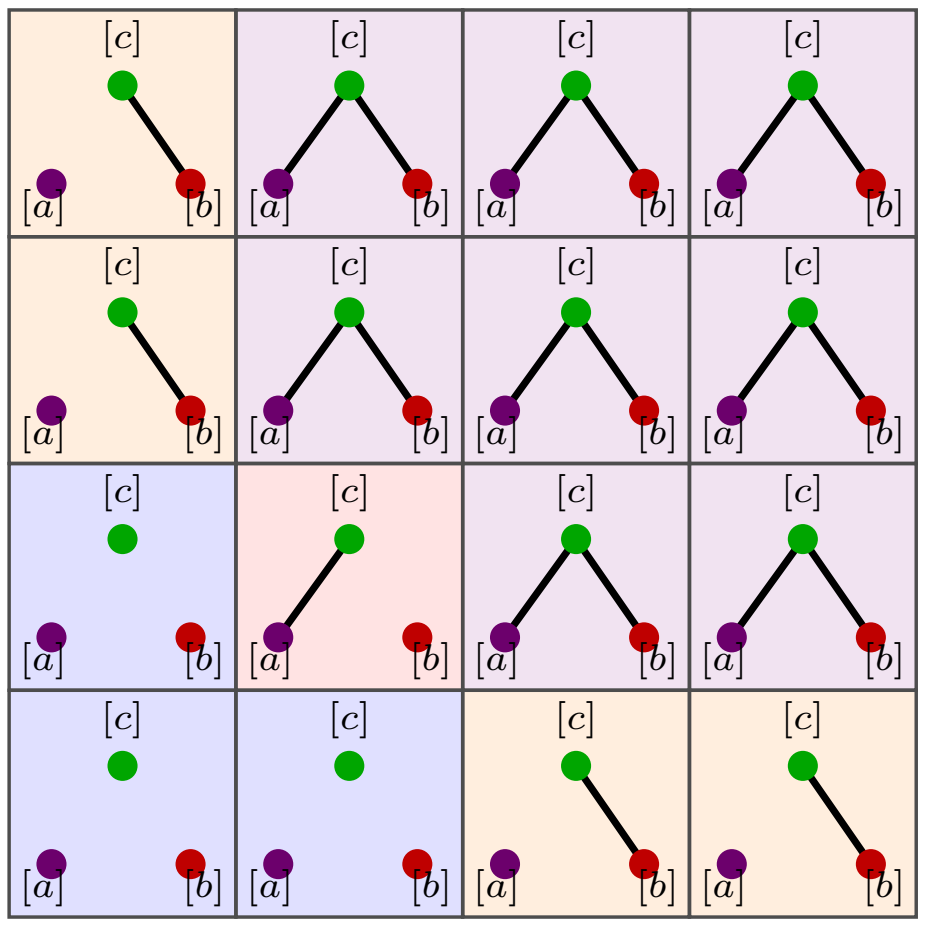
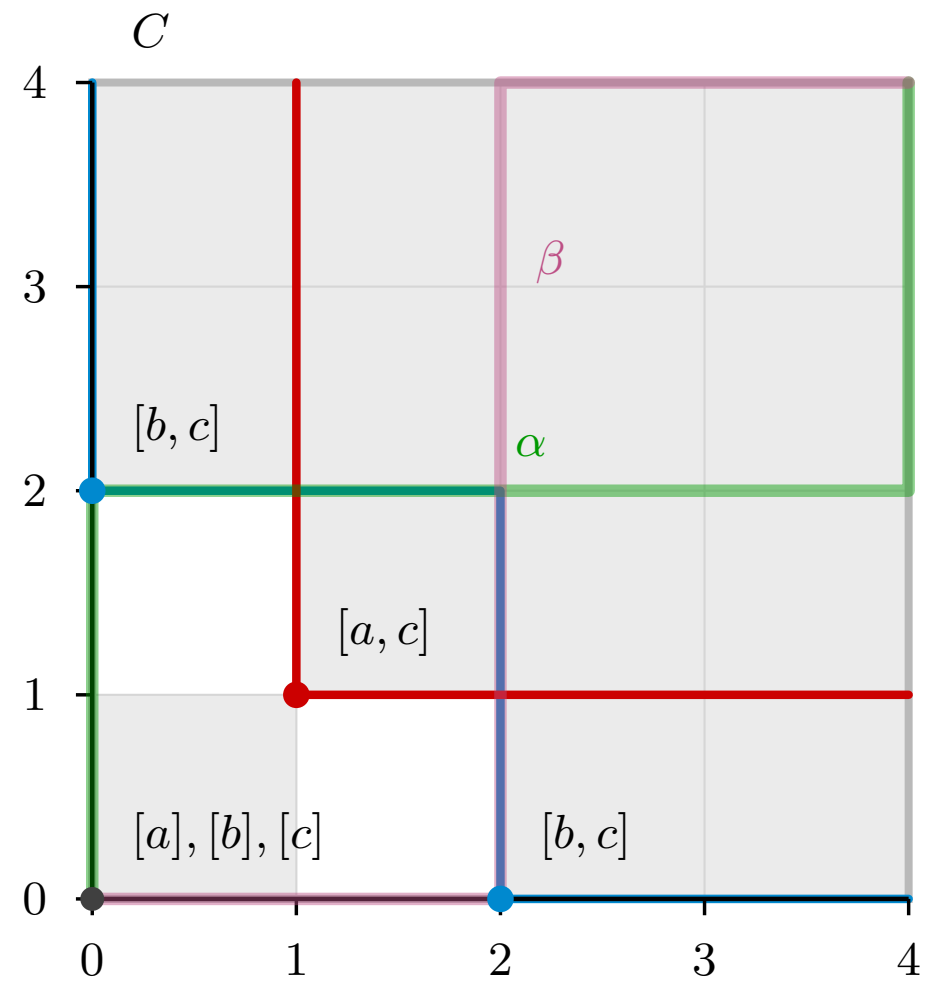
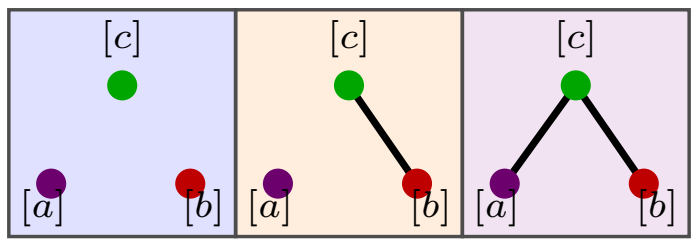
\begin{figure}[h!]
  \centering

  \begin{subfigure}[t]{0.48\textwidth}
    \centering
    \begin{tikzpicture}
      \input{figs/tixzfigures/bifiltration_grid_sketch}
    \end{tikzpicture}
    \caption{A bifiltration $F$.}
  \end{subfigure}
  \hfill
  \begin{subfigure}[t]{0.48\textwidth}
    \centering
    \begin{tikzpicture}
      \input{figs/tixzfigures/directed_parameter_space_complex}
    \end{tikzpicture}
    \caption{Associated entrance complex for $F$, with paths $\alpha$ and $\beta$. The edge $[a,c]$ has entrance coordinate $(1,1)$, and the edge $[b,c]$ has entrance coordinates $(0,2)$ and $(2,0)$.}
  \end{subfigure}

  \medskip

  \begin{subfigure}[t]{0.58\textwidth}
    \centering
    \begin{tikzpicture}
      \input{figs/tixzfigures/directed_alpha_beta_filtrations}
    \end{tikzpicture}
    \caption{The restricted filtration obtained on both $\alpha$ and $\beta$. The paths $\alpha$ and $\beta$ are simplicially equivalent, but not dihomotopy equivalent in $\di{E}$.}
  \end{subfigure}

  \caption{Example filtration with two paths that are simplicially equivalent, but not dihomotopy equivalent in $\di{E}$.}
  \label{fig:se_butnot_de}
\end{figure}

 The converse of \Cref{lem:kcrit_dihomotopy} is not true - it is possible to have simplicially equivalent paths that are not dihomotopy equivalent in $\di{E}$ (see \Cref{fig:se_butnot_de}). However, as simplicial equivalence implies combinatorial equivalence, we have the following corollary.
 
 \begin{corollary} \label{cor:kcrit_dihom}
Let $F$ be a $k$-critical bifiltration with entrance complex $(E,\mathcal{E},\di{E})$, and let $\vecx\prec\vecy\in E$ and let $\alpha,\beta\in \dipaths{E}{\vecx}{\vecy}$.  If there is a dihomotpy $H$ between $\alpha$ and $\beta$ with $H_t\in\dipaths{E}{\vecx}{\vecy}$ for all $t\in[0,1]$, then $\rho(\alpha)=\rho(\beta)$.
 \end{corollary}

\section{Discussion}\label{sec:discussion}

To summarize, in this paper, we have defined a directed space based on the parameter space of a $k$-critical bipersistence module, with the goal of using dihomotopy equivalence to distinguish paths in this complex via combinatorial and/or simplicial equivalence of the barcodes along those paths.  In particular, we have defined the entrance subcomplex (\Cref{def:ec}), a subspace of the directed parameter space complex, in which dihomotopy equivalence of paths ensures both simplicial and combinatorial equivalence along those paths (\Cref{lem:kcrit_dihomotopy}).  Additionally, we have defined an extended fibered barcode (\Cref{def:efb}) for information obtained along non-decreasing paths in parameter space, and shown that the image of the classical fibered barcode is contained within the image of the extended fibered barcode (\Cref{cor:fb_in_efb}).

One important thing to note about this construction is the following - $\di{\mathcal{E}}$ is not a complete representative of $\di{C}$.  In other words, there are paths in $\di{C}$ that do not exist in $\di{\mathcal{E}}$ (such as any path that goes along the diagonal of a mixed cube).  As such, there may be simplicial equivalence classes in the image of the extended fibered barcode that do not have a path representative in $\di{\mathcal{E}}$.  As one of the long-term goals of this directed construction is to compute the extended fibered barcode, and to do so, one must find a representative path for each simplicial equivalence class, this construction may need to be refined in future work.  

Despite this challenge, there are many fruitful directions for future research based on this construction. We list the following conjectures and open questions about our construction:

\begin{conjecture}

It is possible to count the number simplical equivalence classes in the image of the extended fibered barcode by using the information of the entrance complex.

    

\end{conjecture}

In particular, one may obtain a subset of the simplicial equivalence classes in the image of the extended fibered barcode by looking at representative paths in the dihomotopy classes of  $\di{\mathcal{E}}$.  We believe that it is possible to identify a finite subset of paths in $\di C\backslash\di{\mathcal{E}}$ which are representatives for the remaining simplicial equivalence classes by identifying specific paths that pass through mixed cubes, specifically, those which contain the diagonal paths through mixed cubes in $C$.

\begin{conjecture}
It is possible to use directed cubical collapse on $\mathcal{E}$ in order to make computations of $\di{\mathcal{E}}$ (and therefore the extended fibered barcode) more efficient
\end{conjecture}

In past work \cite{belton2020,belton2021combinatorial} we have defined certain sufficient conditions for directed cubical collapse which preserves the topology of the directed path space.  We hope to use these methods on the the entrance complex in a way that simplifies the cubical complex but preserves dihomotopy classes and therefore simplical equivalence classes in the image of the extended fibered barcode.



\section*{Acknowledgements}
This work was partially supported by the Simons Laufer Mathematical Sciences Institute (SLMath) through the PROOF (Promoting Research Opportunities and Open Forums) program, with support from the National Science Foundation and the American Mathematical Society. We really enjoyed our time there, especially the yummy healthy lunches. 

AI disclosure - we used AI for assistance in creating tex for the figures, and to make sure our references are up to date.
\newpage
\bibliography{arxiv_version1}
\bibliographystyle{plainurl}  

\newpage

\end{document}

%% file: figs/tixzfigures/ecc_example.tex
\begin{scope}[
  scale=0.95,
  gridline/.style={draw=gray!35, line width=0.45pt},
  edge/.style={draw=black, line width=1.15pt, line cap=round, line join=round},
  vertex/.style={circle, fill=black, inner sep=1.8pt}
]
  \draw[gridline] (0,0) grid (4,4);

  \fill[gray!45] (0,0) rectangle (1,1);
  \fill[gray!45] (0,1) rectangle (1,2);
  \fill[gray!45] (1,1) rectangle (2,2);
  \fill[gray!45] (2,0) rectangle (3,1);
  \fill[gray!45] (3,3) rectangle (4,4);

  \draw[edge] (0,0) -- (1,0) -- (1,1) -- (2,1) -- (2,2)
              -- (0,2) -- (0,0);
  \draw[edge] (0,1) -- (1,1) -- (1,2);

  \draw[edge] (2,0) rectangle (3,1);

  \draw[edge] (2,2) -- (2,4) -- (4,4) -- (4,3) -- (2,3);
  \draw[edge] (3,3) -- (3,4);

  \foreach \p in {
    (0,0),(1,0),(0,1),(1,1),(2,1),(0,2),(1,2),(2,2),
    (2,0),(3,0),(2,1),(3,1),
    (2,3),(3,3),(4,3),(2,4),(3,4),(4,4)
  } \node[vertex] at \p {};

  \node[anchor=north east, font=\Large] at (0,-0.05) {$0$};
  \node[anchor=north, font=\Large] at (2,-0.08) {$\mathbf{x}$};
  \node[font=\Large] at (2.62,2.60) {$\mathbf{w}$};
  \node[anchor=west, font=\Large] at (4.05,4.05) {$(4,4)$};
\end{scope}

%% file: figs/tixzfigures/bifiltration_grid.tex
\providecommand{\bifiltrationfigurepath}{figs/tixzfigures}

\begin{scope}
\colorlet{emptyHighlight}{black!6}
\colorlet{kzeroHighlight}{blue!18}
\colorlet{koneHighlight}{orange!13}
\colorlet{ktwoHighlight}{cyan!10}
\colorlet{kthreeHighlight}{teal!10}
\colorlet{kfourHighlight}{red!12}
\colorlet{kfiveHighlight}{violet!12}

\newcommand{\complexcell}[2]{%
  \tikz[baseline=(current bounding box.center)]{%
    \fill[#1] (0,0) rectangle (1.5,1.5);%
    \input{\bifiltrationfigurepath/bifiltration_#2}%
  }%
}

\matrix (B) [
  matrix of nodes,
  nodes in empty cells,
  row sep=0pt,
  column sep=0pt,
  nodes={inner sep=0pt, outer sep=0pt}
] {
  %
  \tikz[baseline=(current bounding box.center)]{%
    \fill[kzeroHighlight] (0,0) rectangle (1.5,1.5);%
    \input{\bifiltrationfigurepath/bifiltration_K0}%
  }%
 & %
  \tikz[baseline=(current bounding box.center)]{%
    \fill[kfourHighlight] (0,0) rectangle (1.5,1.5);%
    \input{\bifiltrationfigurepath/bifiltration_K4}%
  }%
 & %
  \tikz[baseline=(current bounding box.center)]{%
    \fill[kfourHighlight] (0,0) rectangle (1.5,1.5);%
    \input{\bifiltrationfigurepath/bifiltration_K4}%
  }%
 & %
  \tikz[baseline=(current bounding box.center)]{%
    \fill[kfourHighlight] (0,0) rectangle (1.5,1.5);%
    \input{\bifiltrationfigurepath/bifiltration_K4}%
  }%
 & %
  \tikz[baseline=(current bounding box.center)]{%
    \fill[kfiveHighlight] (0,0) rectangle (1.5,1.5);%
    \input{\bifiltrationfigurepath/bifiltration_K5}%
  }%
 & %
  \tikz[baseline=(current bounding box.center)]{%
    \fill[kfiveHighlight] (0,0) rectangle (1.5,1.5);%
    \input{\bifiltrationfigurepath/bifiltration_K5}%
  }%
 & %
  \tikz[baseline=(current bounding box.center)]{%
    \fill[kfiveHighlight] (0,0) rectangle (1.5,1.5);%
    \input{\bifiltrationfigurepath/bifiltration_K5}%
  }%
 & %
  \tikz[baseline=(current bounding box.center)]{%
    \fill[kfiveHighlight] (0,0) rectangle (1.5,1.5);%
    \input{\bifiltrationfigurepath/bifiltration_K5}%
  }%
 \\
  %
  \tikz[baseline=(current bounding box.center)]{%
    \fill[kzeroHighlight] (0,0) rectangle (1.5,1.5);%
    \input{\bifiltrationfigurepath/bifiltration_K0}%
  }%
 & %
  \tikz[baseline=(current bounding box.center)]{%
    \fill[kfourHighlight] (0,0) rectangle (1.5,1.5);%
    \input{\bifiltrationfigurepath/bifiltration_K4}%
  }%
 & %
  \tikz[baseline=(current bounding box.center)]{%
    \fill[kfourHighlight] (0,0) rectangle (1.5,1.5);%
    \input{\bifiltrationfigurepath/bifiltration_K4}%
  }%
 & %
  \tikz[baseline=(current bounding box.center)]{%
    \fill[kfourHighlight] (0,0) rectangle (1.5,1.5);%
    \input{\bifiltrationfigurepath/bifiltration_K4}%
  }%
 & %
  \tikz[baseline=(current bounding box.center)]{%
    \fill[kfiveHighlight] (0,0) rectangle (1.5,1.5);%
    \input{\bifiltrationfigurepath/bifiltration_K5}%
  }%
 & %
  \tikz[baseline=(current bounding box.center)]{%
    \fill[kfiveHighlight] (0,0) rectangle (1.5,1.5);%
    \input{\bifiltrationfigurepath/bifiltration_K5}%
  }%
 & %
  \tikz[baseline=(current bounding box.center)]{%
    \fill[kfiveHighlight] (0,0) rectangle (1.5,1.5);%
    \input{\bifiltrationfigurepath/bifiltration_K5}%
  }%
 & %
  \tikz[baseline=(current bounding box.center)]{%
    \fill[kfiveHighlight] (0,0) rectangle (1.5,1.5);%
    \input{\bifiltrationfigurepath/bifiltration_K5}%
  }%
 \\
  %
  \tikz[baseline=(current bounding box.center)]{%
    \fill[kzeroHighlight] (0,0) rectangle (1.5,1.5);%
    \input{\bifiltrationfigurepath/bifiltration_K0}%
  }%
 & %
  \tikz[baseline=(current bounding box.center)]{%
    \fill[kfourHighlight] (0,0) rectangle (1.5,1.5);%
    \input{\bifiltrationfigurepath/bifiltration_K4}%
  }%
 & %
  \tikz[baseline=(current bounding box.center)]{%
    \fill[kfourHighlight] (0,0) rectangle (1.5,1.5);%
    \input{\bifiltrationfigurepath/bifiltration_K4}%
  }%
 & %
  \tikz[baseline=(current bounding box.center)]{%
    \fill[kfourHighlight] (0,0) rectangle (1.5,1.5);%
    \input{\bifiltrationfigurepath/bifiltration_K4}%
  }%
 & %
  \tikz[baseline=(current bounding box.center)]{%
    \fill[kthreeHighlight] (0,0) rectangle (1.5,1.5);%
    \input{\bifiltrationfigurepath/bifiltration_K3}%
  }%
 & %
  \tikz[baseline=(current bounding box.center)]{%
    \fill[kfiveHighlight] (0,0) rectangle (1.5,1.5);%
    \input{\bifiltrationfigurepath/bifiltration_K5}%
  }%
 & %
  \tikz[baseline=(current bounding box.center)]{%
    \fill[kfiveHighlight] (0,0) rectangle (1.5,1.5);%
    \input{\bifiltrationfigurepath/bifiltration_K5}%
  }%
 & %
  \tikz[baseline=(current bounding box.center)]{%
    \fill[kfiveHighlight] (0,0) rectangle (1.5,1.5);%
    \input{\bifiltrationfigurepath/bifiltration_K5}%
  }%
 \\
  %
  \tikz[baseline=(current bounding box.center)]{%
    \fill[kzeroHighlight] (0,0) rectangle (1.5,1.5);%
    \input{\bifiltrationfigurepath/bifiltration_K0}%
  }%
 & %
  \tikz[baseline=(current bounding box.center)]{%
    \fill[kfourHighlight] (0,0) rectangle (1.5,1.5);%
    \input{\bifiltrationfigurepath/bifiltration_K4}%
  }%
 & %
  \tikz[baseline=(current bounding box.center)]{%
    \fill[kfourHighlight] (0,0) rectangle (1.5,1.5);%
    \input{\bifiltrationfigurepath/bifiltration_K4}%
  }%
 & %
  \tikz[baseline=(current bounding box.center)]{%
    \fill[kfourHighlight] (0,0) rectangle (1.5,1.5);%
    \input{\bifiltrationfigurepath/bifiltration_K4}%
  }%
 & %
  \tikz[baseline=(current bounding box.center)]{%
    \fill[kthreeHighlight] (0,0) rectangle (1.5,1.5);%
    \input{\bifiltrationfigurepath/bifiltration_K3}%
  }%
 & %
  \tikz[baseline=(current bounding box.center)]{%
    \fill[kthreeHighlight] (0,0) rectangle (1.5,1.5);%
    \input{\bifiltrationfigurepath/bifiltration_K3}%
  }%
 & %
  \tikz[baseline=(current bounding box.center)]{%
    \fill[kfiveHighlight] (0,0) rectangle (1.5,1.5);%
    \input{\bifiltrationfigurepath/bifiltration_K5}%
  }%
 & %
  \tikz[baseline=(current bounding box.center)]{%
    \fill[kfiveHighlight] (0,0) rectangle (1.5,1.5);%
    \input{\bifiltrationfigurepath/bifiltration_K5}%
  }%
 \\
  %
  \tikz[baseline=(current bounding box.center)]{%
    \fill[kzeroHighlight] (0,0) rectangle (1.5,1.5);%
    \input{\bifiltrationfigurepath/bifiltration_K0}%
  }%
 & %
  \tikz[baseline=(current bounding box.center)]{%
    \fill[koneHighlight] (0,0) rectangle (1.5,1.5);%
    \input{\bifiltrationfigurepath/bifiltration_K1}%
  }%
 & %
  \tikz[baseline=(current bounding box.center)]{%
    \fill[koneHighlight] (0,0) rectangle (1.5,1.5);%
    \input{\bifiltrationfigurepath/bifiltration_K1}%
  }%
 & %
  \tikz[baseline=(current bounding box.center)]{%
    \fill[koneHighlight] (0,0) rectangle (1.5,1.5);%
    \input{\bifiltrationfigurepath/bifiltration_K1}%
  }%
 & %
  \tikz[baseline=(current bounding box.center)]{%
    \fill[ktwoHighlight] (0,0) rectangle (1.5,1.5);%
    \input{\bifiltrationfigurepath/bifiltration_K2}%
  }%
 & %
  \tikz[baseline=(current bounding box.center)]{%
    \fill[ktwoHighlight] (0,0) rectangle (1.5,1.5);%
    \input{\bifiltrationfigurepath/bifiltration_K2}%
  }%
 & %
  \tikz[baseline=(current bounding box.center)]{%
    \fill[ktwoHighlight] (0,0) rectangle (1.5,1.5);%
    \input{\bifiltrationfigurepath/bifiltration_K2}%
  }%
 & %
  \tikz[baseline=(current bounding box.center)]{%
    \fill[ktwoHighlight] (0,0) rectangle (1.5,1.5);%
    \input{\bifiltrationfigurepath/bifiltration_K2}%
  }%
 \\
  %
  \tikz[baseline=(current bounding box.center)]{%
    \fill[kzeroHighlight] (0,0) rectangle (1.5,1.5);%
    \input{\bifiltrationfigurepath/bifiltration_K0}%
  }%
 & %
  \tikz[baseline=(current bounding box.center)]{%
    \fill[koneHighlight] (0,0) rectangle (1.5,1.5);%
    \input{\bifiltrationfigurepath/bifiltration_K1}%
  }%
 & %
  \tikz[baseline=(current bounding box.center)]{%
    \fill[koneHighlight] (0,0) rectangle (1.5,1.5);%
    \input{\bifiltrationfigurepath/bifiltration_K1}%
  }%
 & %
  \tikz[baseline=(current bounding box.center)]{%
    \fill[koneHighlight] (0,0) rectangle (1.5,1.5);%
    \input{\bifiltrationfigurepath/bifiltration_K1}%
  }%
 & %
  \tikz[baseline=(current bounding box.center)]{%
    \fill[ktwoHighlight] (0,0) rectangle (1.5,1.5);%
    \input{\bifiltrationfigurepath/bifiltration_K2}%
  }%
 & %
  \tikz[baseline=(current bounding box.center)]{%
    \fill[ktwoHighlight] (0,0) rectangle (1.5,1.5);%
    \input{\bifiltrationfigurepath/bifiltration_K2}%
  }%
 & %
  \tikz[baseline=(current bounding box.center)]{%
    \fill[ktwoHighlight] (0,0) rectangle (1.5,1.5);%
    \input{\bifiltrationfigurepath/bifiltration_K2}%
  }%
 & %
  \tikz[baseline=(current bounding box.center)]{%
    \fill[ktwoHighlight] (0,0) rectangle (1.5,1.5);%
    \input{\bifiltrationfigurepath/bifiltration_K2}%
  }%
 \\
  %
  \tikz[baseline=(current bounding box.center)]{%
    \fill[kzeroHighlight] (0,0) rectangle (1.5,1.5);%
    \input{\bifiltrationfigurepath/bifiltration_K0}%
  }%
 & %
  \tikz[baseline=(current bounding box.center)]{%
    \fill[kzeroHighlight] (0,0) rectangle (1.5,1.5);%
    \input{\bifiltrationfigurepath/bifiltration_K0}%
  }%
 & %
  \tikz[baseline=(current bounding box.center)]{%
    \fill[koneHighlight] (0,0) rectangle (1.5,1.5);%
    \input{\bifiltrationfigurepath/bifiltration_K1}%
  }%
 & %
  \tikz[baseline=(current bounding box.center)]{%
    \fill[koneHighlight] (0,0) rectangle (1.5,1.5);%
    \input{\bifiltrationfigurepath/bifiltration_K1}%
  }%
 & %
  \tikz[baseline=(current bounding box.center)]{%
    \fill[koneHighlight] (0,0) rectangle (1.5,1.5);%
    \input{\bifiltrationfigurepath/bifiltration_K1}%
  }%
 & %
  \tikz[baseline=(current bounding box.center)]{%
    \fill[koneHighlight] (0,0) rectangle (1.5,1.5);%
    \input{\bifiltrationfigurepath/bifiltration_K1}%
  }%
 & %
  \tikz[baseline=(current bounding box.center)]{%
    \fill[koneHighlight] (0,0) rectangle (1.5,1.5);%
    \input{\bifiltrationfigurepath/bifiltration_K1}%
  }%
 & %
  \tikz[baseline=(current bounding box.center)]{%
    \fill[koneHighlight] (0,0) rectangle (1.5,1.5);%
    \input{\bifiltrationfigurepath/bifiltration_K1}%
  }%
 \\
  %
  \tikz[baseline=(current bounding box.center)]{%
    \fill[emptyHighlight] (0,0) rectangle (1.5,1.5);%
    \input{\bifiltrationfigurepath/bifiltration_Kempty}%
  }%
 & %
  \tikz[baseline=(current bounding box.center)]{%
    \fill[kzeroHighlight] (0,0) rectangle (1.5,1.5);%
    \input{\bifiltrationfigurepath/bifiltration_K0}%
  }%
 & %
  \tikz[baseline=(current bounding box.center)]{%
    \fill[kzeroHighlight] (0,0) rectangle (1.5,1.5);%
    \input{\bifiltrationfigurepath/bifiltration_K0}%
  }%
 & %
  \tikz[baseline=(current bounding box.center)]{%
    \fill[kzeroHighlight] (0,0) rectangle (1.5,1.5);%
    \input{\bifiltrationfigurepath/bifiltration_K0}%
  }%
 & %
  \tikz[baseline=(current bounding box.center)]{%
    \fill[kzeroHighlight] (0,0) rectangle (1.5,1.5);%
    \input{\bifiltrationfigurepath/bifiltration_K0}%
  }%
 & %
  \tikz[baseline=(current bounding box.center)]{%
    \fill[kzeroHighlight] (0,0) rectangle (1.5,1.5);%
    \input{\bifiltrationfigurepath/bifiltration_K0}%
  }%
 & %
  \tikz[baseline=(current bounding box.center)]{%
    \fill[kzeroHighlight] (0,0) rectangle (1.5,1.5);%
    \input{\bifiltrationfigurepath/bifiltration_K0}%
  }%
 & %
  \tikz[baseline=(current bounding box.center)]{%
    \fill[kzeroHighlight] (0,0) rectangle (1.5,1.5);%
    \input{\bifiltrationfigurepath/bifiltration_K0}%
  }%
 \\
};

\foreach \row in {1,...,8} {
  \foreach \col in {1,...,8} {
    \draw[black!70, line width=0.65pt]
      (B-\row-\col.north west) rectangle (B-\row-\col.south east);
  }
}


\end{scope}

%% file: figs/tixzfigures/biparameterspace_example.tex
\begin{scope}[x=1.15cm,y=1.15cm]
  \colorlet{supportA}{blue!80!black}
  \colorlet{supportBCA}{orange!90!black}
  \colorlet{supportAC}{red!80!black}
  \colorlet{supportBC}{cyan!70!blue}
  \colorlet{supportABC}{violet!85!black}

  \tikzset{
    support boundary/.style={line cap=round,line join=round},
    entrance point/.style={circle,inner sep=0pt}
  }

  \draw[support boundary, supportA, line width=3pt]
    (0,7) -- (0,1) -- (1,1) -- (1,0) -- (7,0);
  \draw[black, line width=.85pt] (0,7) -- (0,0) -- (7,0);
  \foreach \t in {0,...,7} {
    \draw[black, line width=.65pt] (\t,0) -- (\t,-.08);
    \draw[black, line width=.65pt] (0,\t) -- (-.08,\t);
  }

  \node[anchor=north, font=\small] at (0,-0.13) {$0$};
  \node[anchor=north, font=\small] at (1,-0.13) {$1$};
  \node[anchor=north, font=\small] at (2,-0.13) {$2$};
  \node[anchor=north, font=\small] at (3,-0.13) {$3$};
  \node[anchor=north, font=\small] at (4,-0.13) {$4$};
  \node[anchor=north, font=\small] at (5,-0.13) {$5$};
  \node[anchor=north, font=\small] at (6,-0.13) {$6$};
  \node[anchor=north, font=\small] at (7,-0.13) {$7$};
  \node[anchor=east, font=\small] at (-0.13,0) {$0$};
  \node[anchor=east, font=\small] at (-0.13,1) {$1$};
  \node[anchor=east, font=\small] at (-0.13,2) {$2$};
  \node[anchor=east, font=\small] at (-0.13,3) {$3$};
  \node[anchor=east, font=\small] at (-0.13,4) {$4$};
  \node[anchor=east, font=\small] at (-0.13,5) {$5$};
  \node[anchor=east, font=\small] at (-0.13,6) {$6$};
  \node[anchor=east, font=\small] at (-0.13,7) {$7$};

  \draw[support boundary, supportBCA, line width=2.2pt]
    (1,7) -- (1,2) -- (2,2) -- (2,1) -- (7,1);

  \draw[support boundary, supportBC, line width=2.2pt]
    (4,7) -- (4,2) -- (7,2);
  \draw[support boundary, supportAC, line width=2.2pt]
    (1,7) -- (1,4) -- (7,4);
  \draw[support boundary, supportABC, line width=2.2pt]
    (4,7) -- (4,6) -- (5,6) -- (5,5) -- (6,5) -- (6,4) -- (7,4);

  \foreach \p in {(0,1),(1,0)} {
    \node[entrance point, fill=supportA, minimum size=6pt] at \p {};
  }

  \foreach \p in {(1,2),(2,1)} {
    \node[entrance point, fill=supportBCA, minimum size=6pt] at \p {};
  }

  \node[entrance point, fill=supportBC, minimum size=6pt] at (4,2) {};
  \node[entrance point, fill=supportAC, minimum size=6pt] at (1,4) {};
  \foreach \p in {(4,6),(5,5),(6,4)} {
    \node[entrance point, fill=supportABC, minimum size=6pt] at \p {};
  }

  \node[anchor=south west, font=\small] at (0.10,1.10) {$[a]$};
  \node[anchor=south west, font=\small] at (1.10,2.10) {$[b],[c],[a,b]$};
  \node[anchor=south west, font=\small] at (4.10,2.10) {$[a,c]$};
  \node[anchor=south west, font=\small] at (1.10,4.10) {$[b,c]$};
  \node[anchor=south west, font=\small] at (5.10,5.10) {$[a,b,c]$};
\end{scope}

%% file: figs/tixzfigures/directedPS_example.tex
\begin{scope}[x=1.15cm,y=1.15cm]
  \colorlet{supportA}{blue!80!black}
  \colorlet{supportBCA}{orange!90!black}
  \colorlet{supportAC}{red!80!black}
  \colorlet{supportBC}{cyan!70!blue}
  \colorlet{supportABC}{violet!85!black}

  \tikzset{
    support boundary/.style={line cap=round,line join=round},
    entrance point/.style={circle,inner sep=0pt}
  }

  \draw[gray!55, line width=1.5pt] (0,0) rectangle (6,6);
 \foreach \t in {1,...,7} {
    \draw[black!16, line width=.4pt] (\t,0) -- (\t,7);
    \draw[black!16, line width=.4pt] (0,\t) -- (7,\t);
  }
  \draw[support boundary, supportA, line width=3pt]
    (0,7) -- (0,1) -- (1,1) -- (1,0) -- (7,0);
  \draw[black, line width=.85pt] (0,7) -- (0,0) -- (7,0);
  \foreach \t in {0,...,7} {
    \draw[black, line width=.65pt] (\t,0) -- (\t,-.08);
    \draw[black, line width=.65pt] (0,\t) -- (-.08,\t);
  }

  \node[anchor=north, font=\small] at (0,-0.13) {$0$};
  \node[anchor=north, font=\small] at (1,-0.13) {$1$};
  \node[anchor=north, font=\small] at (2,-0.13) {$2$};
  \node[anchor=north, font=\small] at (3,-0.13) {$3$};
  \node[anchor=north, font=\small] at (4,-0.13) {$4$};
  \node[anchor=north, font=\small] at (5,-0.13) {$5$};
  \node[anchor=north, font=\small] at (6,-0.13) {$6$};
  \node[anchor=north, font=\small] at (7,-0.13) {$7$};
  \node[anchor=east, font=\small] at (-0.13,0) {$0$};
  \node[anchor=east, font=\small] at (-0.13,1) {$1$};
  \node[anchor=east, font=\small] at (-0.13,2) {$2$};
  \node[anchor=east, font=\small] at (-0.13,3) {$3$};
  \node[anchor=east, font=\small] at (-0.13,4) {$4$};
  \node[anchor=east, font=\small] at (-0.13,5) {$5$};
  \node[anchor=east, font=\small] at (-0.13,6) {$6$};
  \node[anchor=east, font=\small] at (-0.13,7) {$7$};

  \draw[support boundary, supportBCA, line width=2.2pt]
    (1,7) -- (1,2) -- (2,2) -- (2,1) -- (7,1);

  \draw[support boundary, supportBC, line width=2.2pt]
    (4,7) -- (4,2) -- (7,2);
  \draw[support boundary, supportAC, line width=2.2pt]
    (1,7) -- (1,4) -- (7,4);
  \draw[support boundary, supportABC, line width=2.2pt]
    (4,7) -- (4,6) -- (5,6) -- (5,5) -- (6,5) -- (6,4) -- (7,4);

  \foreach \p in {(0,1),(1,0)} {
    \node[entrance point, fill=supportA, minimum size=6pt] at \p {};
  }

  \foreach \p in {(1,2),(2,1)} {
    \node[entrance point, fill=supportBCA, minimum size=6pt] at \p {};
  }

  \node[entrance point, fill=supportBC, minimum size=6pt] at (4,2) {};
  \node[entrance point, fill=supportAC, minimum size=6pt] at (1,4) {};
  \foreach \p in {(4,6),(5,5),(6,4)} {
    \node[entrance point, fill=supportABC, minimum size=6pt] at \p {};
  }

  \node[anchor=south west, font=\small] at (0.10,1.10) {$[a]$};
  \node[anchor=south west, font=\small] at (1.10,2.10) {$[b],[c],[a,b]$};
  \node[anchor=south west, font=\small] at (4.10,2.10) {$[a,c]$};
  \node[anchor=south west, font=\small] at (1.10,4.10) {$[b,c]$};
  \node[anchor=south west, font=\small] at (5.10,5.10) {$[a,b,c]$};
\end{scope}

%% file: figs/tixzfigures/path_in_dpsc_ex.tex
\begin{scope}[x=1.15cm,y=1.15cm]
  \colorlet{supportA}{blue!80!black}
  \colorlet{supportBCA}{orange!90!black}
  \colorlet{supportAC}{red!80!black}
  \colorlet{supportBC}{cyan!70!blue}
  \colorlet{supportABC}{violet!85!black}
  \colorlet{pathAlpha}{green!60!black}

  \tikzset{
    support boundary/.style={line cap=round,line join=round},
    entrance point/.style={circle,inner sep=0pt}
  }

  \foreach \t in {1,...,7} {
    \draw[black!16, line width=.4pt] (\t,0) -- (\t,7);
    \draw[black!16, line width=.4pt] (0,\t) -- (7,\t);
  }

  \draw[gray!55, line width=1.5pt] (0,0) rectangle (6,6);

  \draw[support boundary, supportA, line width=3pt]
    (0,7) -- (0,1) -- (1,1) -- (1,0) -- (7,0);

  \draw[black, line width=.85pt] (0,7) -- (0,0) -- (7,0);
  \foreach \t in {0,...,7} {
    \draw[black, line width=.65pt] (\t,0) -- (\t,-.08);
    \draw[black, line width=.65pt] (0,\t) -- (-.08,\t);
  }

  \node[anchor=north, font=\small] at (0,-0.13) {$0$};
  \node[anchor=north, font=\small] at (1,-0.13) {$1$};
  \node[anchor=north, font=\small] at (2,-0.13) {$2$};
  \node[anchor=north, font=\small] at (3,-0.13) {$3$};
  \node[anchor=north, font=\small] at (4,-0.13) {$4$};
  \node[anchor=north, font=\small] at (5,-0.13) {$5$};
  \node[anchor=north, font=\small] at (6,-0.13) {$6$};
  \node[anchor=north, font=\small] at (7,-0.13) {$7$};
  \node[anchor=east, font=\small] at (-0.13,0) {$0$};
  \node[anchor=east, font=\small] at (-0.13,1) {$1$};
  \node[anchor=east, font=\small] at (-0.13,2) {$2$};
  \node[anchor=east, font=\small] at (-0.13,3) {$3$};
  \node[anchor=east, font=\small] at (-0.13,4) {$4$};
  \node[anchor=east, font=\small] at (-0.13,5) {$5$};
  \node[anchor=east, font=\small] at (-0.13,6) {$6$};
  \node[anchor=east, font=\small] at (-0.13,7) {$7$};

  \draw[support boundary, supportBCA, line width=1.6pt]
    (1,7) -- (1,2) -- (2,2) -- (2,1) -- (7,1);

  \draw[support boundary, supportBC, line width=1.6pt]
    (4,7) -- (4,2) -- (7,2);
  \draw[support boundary, supportAC, line width=1.6pt]
    (1,7) -- (1,4) -- (7,4);
  \draw[support boundary, supportABC, line width=1.6pt]
    (4,7) -- (4,6) -- (5,6) -- (5,5) -- (6,5) -- (6,4) -- (7,4);

  \draw[pathAlpha, line width=2.35pt, line cap=round, line join=round]
    (0,0)
    .. controls (0.22,0.08) and (1.10,1.82) .. (1.45,2)
    .. controls (1.80,2.18) and (2.30,3.70) .. (2.65,4)
    .. controls (3.00,4.30) and (3.60,4.40) .. (4,4.45)
    .. controls (4.40,4.50) and (4.60,5.05) .. (5,5.25)
    .. controls (5.40,5.45) and (5.72,5.86) .. (6,6);

  \node[fill=white, inner sep=1pt, font=\small, text=pathAlpha]
    at (3.35,4.72) {$\alpha$};

  \node[entrance point, fill=supportA, minimum size=6pt] at (0.6996,1) {};
  \node[entrance point, fill=supportBCA, minimum size=6pt] at (1.45,2) {};
  \node[entrance point, fill=supportAC, minimum size=6pt] at (2.65,4) {};
  \node[entrance point, fill=supportBC, minimum size=6pt] at (4,4.45) {};
  \node[entrance point, fill=supportABC, minimum size=6pt] at (5,5.25) {};

  \node[anchor=south east, font=\small, text=supportAC]
    at (3,4.2) {$\push_{\alpha}([b,c])$};

  \foreach \p in {(0,1),(1,0)} {
    \node[entrance point, fill=supportA, minimum size=6pt] at \p {};
  }

  \foreach \p in {(1,2),(2,1)} {
    \node[entrance point, fill=supportBCA, minimum size=6pt] at \p {};
  }

  \node[entrance point, fill=supportBC, minimum size=6pt] at (4,2) {};
  \node[entrance point, fill=supportAC, minimum size=6pt] at (1,4) {};
  \foreach \p in {(4,6),(5,5),(6,4)} {
    \node[entrance point, fill=supportABC, minimum size=6pt] at \p {};
  }

\end{scope}

%% file: figs/tixzfigures/path_dpsc_filtration.tex
\providecommand{\bifiltrationfigurepath}{figs/tixzfigures}

\begin{scope}
  \newcommand{\alphacell}[1]{%
    \tikz[baseline=(current bounding box.center)]{%
      \input{\bifiltrationfigurepath/bifiltration_#1}%
    }%
  }

  \matrix (A) [
    matrix of nodes,
    nodes in empty cells,
    row sep=0pt,
    column sep=0pt,
    nodes={inner sep=0pt, outer sep=0pt}
  ] {
    \tikz[baseline=(current bounding box.center)]{%
      \input{\bifiltrationfigurepath/bifiltration_Kempty}%
    }%
  &
    \tikz[baseline=(current bounding box.center)]{%
      \input{\bifiltrationfigurepath/bifiltration_K0}%
    }%
   &
    \tikz[baseline=(current bounding box.center)]{%
      \input{\bifiltrationfigurepath/bifiltration_K1}%
    }%
   &
    \tikz[baseline=(current bounding box.center)]{%
      \input{\bifiltrationfigurepath/bifiltration_K4}%
    }%
   &
    \tikz[baseline=(current bounding box.center)]{%
      \input{\bifiltrationfigurepath/bifiltration_K3}%
    }%
   &
    \tikz[baseline=(current bounding box.center)]{%
      \input{\bifiltrationfigurepath/bifiltration_K5}%
    }%
   \\
  };

  \foreach \col in {1,...,6} {
    \draw[black!70, line width=0.65pt]
      (A-1-\col.north west) rectangle (A-1-\col.south east);
  }
\end{scope}

%% file: figs/tixzfigures/barcode_A.tex
\begin{scope}[x=.43cm,y=.68cm]
  \colorlet{barcodeHzero}{violet!80!black}
  \colorlet{barcodeHone}{blue!70!black}

  \draw[black!70, line width=.55pt, ->] (0,.15) -- (11,.15);
  \draw[black!10, dashed, line width=.5pt] (0,1.72) -- (11,1.72);
  \foreach \t in {0,...,10} {
    \draw[black!60, line width=.45pt] (\t,.15) -- (\t,.07);
    \node[anchor=north, font=\scriptsize, text=black!70] at (\t,.02) {$\t$};
  }

  \node[anchor=east, font=\small] at (-.55,.78) {$H_0$};
  \node[anchor=east, font=\small] at (-.55,2.62) {$H_1$};

  \draw[barcodeHzero, line width=1.35pt, ->] (0,.78) -- (10.25,.78);
  \fill[barcodeHzero] (0,.78) circle (2.5pt);

  \draw[barcodeHone, line width=1.35pt] (1,2.25) -- (5,2.25);
  \fill[barcodeHone] (1,2.25) circle (2.5pt);
  \fill[barcodeHone] (5,2.25) circle (2.5pt);
  \draw[barcodeHone, line width=1.35pt] (6,2.85) -- (8,2.85);
  \fill[barcodeHone] (6,2.85) circle (2.5pt);
  \fill[barcodeHone] (8,2.85) circle (2.5pt);

\end{scope}

%% file: figs/tixzfigures/barcode_rhoA.tex
\begin{scope}[x=.43cm,y=.68cm]
  \colorlet{barcodeHzero}{violet!80!black}
  \colorlet{barcodeHone}{blue!70!black}

  \draw[black!70, line width=.55pt, ->] (0,.15) -- (11,.15);
  \draw[black!10, dashed, line width=.5pt] (0,1.72) -- (11,1.72);
  \foreach \t in {0,...,10} {
    \draw[black!60, line width=.45pt] (\t,.15) -- (\t,.07);
    \node[anchor=north, font=\scriptsize, text=black!70] at (\t,.02) {$\t$};
  }

  \node[anchor=east, font=\small] at (-.55,.78) {$H_0$};
  \node[anchor=east, font=\small] at (-.55,2.62) {$H_1$};

  \draw[barcodeHzero, line width=1.35pt, ->] (0,.78) -- (10.25,.78);
  \fill[barcodeHzero] (0,.78) circle (2.5pt);

  \draw[barcodeHone, line width=1.35pt] (1,2.25) -- (2,2.25);
  \fill[barcodeHone] (1,2.25) circle (2.5pt);
  \fill[barcodeHone] (2,2.25) circle (2.5pt);
  \draw[barcodeHone, line width=1.35pt] (3,2.85) -- (4,2.85);
  \fill[barcodeHone] (3,2.85) circle (2.5pt);
  \fill[barcodeHone] (4,2.85) circle (2.5pt);
\end{scope}

%% file: figs/tixzfigures/barcode_B.tex
\begin{scope}[x=.43cm,y=.68cm]
  \colorlet{barcodeHzero}{violet!80!black}
  \colorlet{barcodeHone}{blue!70!black}

  \draw[black!70, line width=.55pt, ->] (0,.15) -- (11,.15);
  \draw[black!10, dashed, line width=.5pt] (0,1.72) -- (11,1.72);
  \foreach \t in {0,...,10} {
    \draw[black!60, line width=.45pt] (\t,.15) -- (\t,.07);
    \node[anchor=north, font=\scriptsize, text=black!70] at (\t,.02) {$\t$};
  }

  \node[anchor=east, font=\small] at (-.55,.78) {$H_0$};
  \node[anchor=east, font=\small] at (-.55,2.62) {$H_1$};

  \draw[barcodeHzero, line width=1.35pt, ->] (0,.78) -- (10.25,.78);
  \fill[barcodeHzero] (0,.78) circle (2.5pt);

  \draw[barcodeHone, line width=1.35pt] (1,2.25) -- (5.2,2.25);
  \fill[barcodeHone] (1,2.25) circle (2.5pt);
  \fill[barcodeHone] (5.2,2.25) circle (2.5pt);
  \draw[barcodeHone, line width=1.35pt] (7.2,2.85) -- (8.2,2.85);
  \fill[barcodeHone] (7.2,2.85) circle (2.5pt);
  \fill[barcodeHone] (8.2,2.85) circle (2.5pt);

\end{scope}

%% file: figs/tixzfigures/barcode_D.tex
\begin{scope}[x=.43cm,y=.68cm]
  \colorlet{barcodeHzero}{violet!80!black}
  \colorlet{barcodeHone}{blue!70!black}

  \draw[black!70, line width=.55pt, ->] (0,.15) -- (11,.15);
  \draw[black!10, dashed, line width=.5pt] (0,1.72) -- (11,1.72);
  \foreach \t in {0,...,10} {
    \draw[black!60, line width=.45pt] (\t,.15) -- (\t,.07);
    \node[anchor=north, font=\scriptsize, text=black!70] at (\t,.02) {$\t$};
  }

  \node[anchor=east, font=\small] at (-.55,.78) {$H_0$};
  \node[anchor=east, font=\small] at (-.55,2.62) {$H_1$};

  \draw[barcodeHzero, line width=1.35pt, ->] (0,.78) -- (10.25,.78);
  \fill[barcodeHzero] (0,.78) circle (2.5pt);

  \draw[barcodeHone, line width=1.35pt] (1,2.25) -- (5,2.25);
  \fill[barcodeHone] (1,2.25) circle (2.5pt);
  \fill[barcodeHone] (5,2.25) circle (2.5pt);
  \draw[barcodeHone, line width=1.35pt] (5,2.85) -- (8,2.85);
  \fill[barcodeHone] (5,2.85) circle (2.5pt);
  \fill[barcodeHone] (8,2.85) circle (2.5pt);
\end{scope}

%% file: figs/tixzfigures/incomparablelemma_ex.tex
\begin{scope}[x=1.15cm,y=1.15cm]
  \colorlet{supportA}{blue!80!black}
  \colorlet{supportBCA}{orange!90!black}
  \colorlet{supportAC}{red!80!black}
  \colorlet{supportBC}{cyan!70!blue}
  \colorlet{supportABC}{violet!85!black}

  \colorlet{pathAlpha}{green!60!black}
  \colorlet{pathBeta}{black!75}
  \colorlet{pathGamma}{magenta!75!black}

  \tikzset{
    support boundary/.style={line cap=round,line join=round},
    entrance point/.style={circle,inner sep=0pt},
    filtration path/.style={line width=2.35pt,line cap=round,line join=round}
  }

  \foreach \t in {1,...,7} {
    \draw[black!16, line width=.4pt] (\t,0) -- (\t,7);
    \draw[black!16, line width=.4pt] (0,\t) -- (7,\t);
  }

  \draw[gray!55, line width=1.5pt] (0,0) rectangle (6,6);

  \draw[support boundary, supportA, line width=3pt]
    (0,7) -- (0,1) -- (1,1) -- (1,0) -- (7,0);

  \draw[black, line width=.85pt] (0,7) -- (0,0) -- (7,0);
  \foreach \t in {0,...,7} {
    \draw[black, line width=.65pt] (\t,0) -- (\t,-.08);
    \draw[black, line width=.65pt] (0,\t) -- (-.08,\t);
  }

  \node[anchor=north, font=\small] at (0,-0.13) {$0$};
  \node[anchor=north, font=\small] at (1,-0.13) {$1$};
  \node[anchor=north, font=\small] at (2,-0.13) {$2$};
  \node[anchor=north, font=\small] at (3,-0.13) {$3$};
  \node[anchor=north, font=\small] at (4,-0.13) {$4$};
  \node[anchor=north, font=\small] at (5,-0.13) {$5$};
  \node[anchor=north, font=\small] at (6,-0.13) {$6$};
  \node[anchor=north, font=\small] at (7,-0.13) {$7$};
  \node[anchor=east, font=\small] at (-0.13,0) {$0$};
  \node[anchor=east, font=\small] at (-0.13,1) {$1$};
  \node[anchor=east, font=\small] at (-0.13,2) {$2$};
  \node[anchor=east, font=\small] at (-0.13,3) {$3$};
  \node[anchor=east, font=\small] at (-0.13,4) {$4$};
  \node[anchor=east, font=\small] at (-0.13,5) {$5$};
  \node[anchor=east, font=\small] at (-0.13,6) {$6$};
  \node[anchor=east, font=\small] at (-0.13,7) {$7$};

  \draw[support boundary, supportBCA, line width=1.6pt]
    (1,7) -- (1,2) -- (2,2) -- (2,1) -- (7,1);

  \draw[support boundary, supportBC, line width=1.6pt]
    (4,7) -- (4,2) -- (7,2);
  \draw[support boundary, supportAC, line width=1.6pt]
    (1,7) -- (1,4) -- (7,4);
  \draw[support boundary, supportABC, line width=1.6pt]
    (4,7) -- (4,6) -- (5,6) -- (5,5) -- (6,5) -- (6,4) -- (7,4);

  \draw[filtration path, pathAlpha]
    (0,0)
    .. controls (0.22,0.08) and (1.10,1.82) .. (1.45,2)
    .. controls (1.80,2.18) and (2.30,3.70) .. (2.65,4)
    .. controls (3.00,4.30) and (3.60,4.40) .. (4,4.45)
    .. controls (4.40,4.50) and (4.60,5.05) .. (5,5.25)
    .. controls (5.40,5.45) and (5.72,5.86) .. (6,6);
  \node[font=\small, text=pathAlpha] at (3.35,3.72) {$\alpha$};

  \draw[filtration path, pathBeta]
    (0,0)
    .. controls (0.17,0.13) and (0.47,0.70) .. (0.55,1)
    .. controls (0.63,1.30) and (0.97,4.12) .. (1,4.42)
    .. controls (1.03,4.52) and (3.73,4.76) .. (4,4.82)
    .. controls (4.27,4.88) and (4.77,5.11) .. (5,5.28)
    .. controls (5.23,5.45) and (5.83,5.83) .. (6,6);
  \node[font=\small, text=pathBeta] at (2.25,5.10) {$\beta$};

  \draw[filtration path, pathGamma]
    (0,0)
    .. controls (0.17,0.13) and (0.35,0.83) .. (0.55,1)
    .. controls (0.75,1.17) and (1.30,1.83) .. (1.60,2)
    .. controls (1.90,2.10) and (3.83,2.40) .. (4,2.50)
    .. controls (4.17,2.60) and (4.42,3.73) .. (4.50,4)
    .. controls (4.58,4.27) and (4.80,5.20) .. (5,5.40)
    .. controls (5.20,5.60) and (5.83,5.87) .. (6,6);
  \node[font=\small, text=pathGamma] at (2.70,2.85) {$\gamma$};
  
  \foreach \p in {(0,1),(1,0)} {
    \node[entrance point, fill=supportA, minimum size=6pt] at \p {};
  }

  \foreach \p in {(1,2),(2,1)} {
    \node[entrance point, fill=supportBCA, minimum size=6pt] at \p {};
  }

  \node[entrance point, fill=supportBC, minimum size=6pt] at (4,2) {};
  \node[entrance point, fill=supportAC, minimum size=6pt] at (1,4) {};
  \foreach \p in {(4,6),(5,5),(6,4)} {
    \node[entrance point, fill=supportABC, minimum size=6pt] at \p {};
  }

\end{scope}

%% file: figs/tixzfigures/efb_example.tex
\begin{scope}[x=.72cm,y=.72cm]
  \colorlet{supportA}{blue!80!black}
  \colorlet{supportBCA}{orange!90!black}
  \colorlet{supportAC}{red!80!black}
  \colorlet{supportBC}{cyan!70!blue}
  \colorlet{supportABC}{violet!85!black}

  \colorlet{lineOne}{black!75}
  \colorlet{lineTwo}{green!60!black}
  \colorlet{lineThree}{magenta!75!black}

  \tikzset{
    support boundary/.style={line cap=round,line join=round},
    entrance point/.style={circle,inner sep=0pt},
    affine line/.style={line width=2pt,line cap=round},
 highlighted path/.style={
  line width=3pt,
  line cap=round,
  line join=round,
  opacity=.3
}
  }

  \foreach \t in {1,...,12} {
    \draw[black!16, line width=.4pt] (\t,0) -- (\t,12);
    \draw[black!16, line width=.4pt] (0,\t) -- (12,\t);
  }

  \draw[gray!55, line width=1.5pt] (0,0) rectangle (6,6);
  \node[anchor=south west, font=\small] at (0.12,6.08) {$C$};

  \draw[support boundary, supportA, line width=1.55pt]
    (0,12) -- (0,1) -- (1,1) -- (1,0) -- (12,0);

  \draw[black, line width=.85pt] (0,12) -- (0,0) -- (12,0);
  \foreach \t in {0,...,12} {
    \draw[black, line width=.65pt] (\t,0) -- (\t,-.10);
    \draw[black, line width=.65pt] (0,\t) -- (-.10,\t);
  }

  \foreach \t in {0,...,12} {
    \node[anchor=north, font=\small] at (\t,-0.16) {$\t$};
    \node[anchor=east, font=\small] at (-0.16,\t) {$\t$};
  }

  \draw[support boundary, supportBCA, line width=1.1pt]
    (1,12) -- (1,2) -- (2,2) -- (2,1) -- (12,1);

  \draw[support boundary, supportBC, line width=1.1pt]
    (4,12) -- (4,2) -- (12,2);
  \draw[support boundary, supportAC, line width=1.1pt]
    (1,12) -- (1,4) -- (12,4);
  \draw[support boundary, supportABC, line width=1.1pt]
    (4,12) -- (4,6) -- (5,6) -- (5,5) -- (6,5) -- (6,4) -- (12,4);

  \draw[affine line, lineOne]
    (0,2) -- (12,8)
    node[pos=.76, above, font=\small, text=lineOne] {$L_1$};

  \draw[affine line, lineTwo]
    (0,6.25) -- (12,11)
    node[pos=.74, above, font=\small, text=lineTwo] {$L_2$};

  \draw[affine line, lineThree]
    (8,0) -- (12,7)
    node[pos=.52, below, font=\small, text=lineThree] {$L_3$};

  \draw[highlighted path, lineOne]
    (0,2) -- (6,5) -- (6,6);
  \node[font=\small, text=lineOne] at (2.55,3.75) {$\alpha_{L_1}$};

  \draw[highlighted path, lineThree] (6,0) -- (6,6);
  \node[anchor=west, font=\small, text=lineThree] at (6,3) {$\alpha_{L_3}$};

  \draw[highlighted path, lineTwo] (0,6) -- (6,6);
  \node[anchor=south, font=\small, text=lineTwo] at (3,6) {$\alpha_{L_2}$};

  \foreach \p in {(0,1),(1,0)} {
    \node[entrance point, fill=supportA, minimum size=4.2pt] at \p {};
  }

  \foreach \p in {(1,2),(2,1)} {
    \node[entrance point, fill=supportBCA, minimum size=4.2pt] at \p {};
  }

  \node[entrance point, fill=supportBC, minimum size=4.2pt] at (4,2) {};
  \node[entrance point, fill=supportAC, minimum size=4.2pt] at (1,4) {};
  \foreach \p in {(4,6),(5,5),(6,4)} {
    \node[entrance point, fill=supportABC, minimum size=4.2pt] at \p {};
  }

\end{scope}

%% file: figs/tixzfigures/six_thresholds_alpha.tex
\begin{scope}[x=1.0cm,y=1.0cm]
  \colorlet{thresholdOne}{blue!80!black}
  \colorlet{thresholdTwo}{orange!90!black}
  \colorlet{thresholdThree}{red!80!black}
  \colorlet{thresholdFour}{cyan!70!blue}
  \colorlet{thresholdFive}{violet!85!black}
  \colorlet{thresholdSix}{green!60!black}
  \colorlet{lineL}{gray!70}
  \colorlet{pathAlpha}{black!85}

  \tikzset{
    support boundary/.style={line cap=round,line join=round},
    entrance point/.style={circle,inner sep=0pt}
  }

  \foreach \t in {1,...,6} {
    \draw[black!16, line width=.4pt] (\t,0) -- (\t,6);
    \draw[black!16, line width=.4pt] (0,\t) -- (6,\t);
  }

  \draw[gray!55, line width=1.5pt] (0,0) rectangle (6,6);
  \node[anchor=south west, font=\small] at (0.10,6.08) {$C$};

  \draw[support boundary, thresholdOne, line width=1.55pt]
    (0,6) -- (0,1) -- (6,1);
  \draw[support boundary, thresholdTwo, line width=1.55pt]
    (1,6) -- (1,0) -- (6,0);
  \draw[support boundary, thresholdThree, line width=1.55pt]
    (1,6) -- (1,2) -- (6,2);
  \draw[support boundary, thresholdFour, line width=1.55pt]
    (2,6) -- (2,1) -- (6,1);
  \draw[support boundary, thresholdFive, line width=1.55pt]
    (3,6) -- (3,2) -- (6,2);
  \draw[support boundary, thresholdSix, line width=1.55pt]
    (2,6) -- (2,3) -- (6,3);

  \draw[black, line width=.85pt] (0,6) -- (0,0) -- (6,0);
  \foreach \t in {0,...,6} {
    \draw[black, line width=.65pt] (\t,0) -- (\t,-.08);
    \draw[black, line width=.65pt] (0,\t) -- (-.08,\t);
  }

  \foreach \t in {0,...,6} {
    \node[anchor=north, font=\small] at (\t,-0.13) {$\t$};
    \node[anchor=east, font=\small] at (-0.13,\t) {$\t$};
  }

  \draw[lineL, line width=2.2pt, line cap=round] (0,0) -- (6,6);
  \node[anchor=north west, font=\small, text=lineL] at (4.05,4.00) {$L$};

  \draw[pathAlpha, line width=2.5pt, line cap=round, line join=round]
    (0,0) -- (2,2) -- (2,3) -- (2,6) -- (6,6);
  \node[anchor=west, font=\small, text=pathAlpha] at (2.10,4.55) {$\alpha$};

  \node[entrance point, fill=thresholdOne, minimum size=5pt] at (0,1) {};
  \node[entrance point, fill=thresholdTwo, minimum size=5pt] at (1,0) {};
  \node[entrance point, fill=thresholdThree, minimum size=5pt] at (1,2) {};
  \node[entrance point, fill=thresholdFour, minimum size=5pt] at (2,1) {};
  \node[entrance point, fill=thresholdFive, minimum size=5pt] at (3,2) {};
  \node[entrance point, fill=thresholdSix, minimum size=5pt] at (2,3) {};
\end{scope}

%% file: figs/tixzfigures/entrance_subcomplex.tex
\begin{scope}[x=1.15cm,y=1.15cm]
  \colorlet{supportA}{blue!80!black}
  \colorlet{supportBCA}{orange!90!black}
  \colorlet{supportAC}{red!80!black}
  \colorlet{supportBC}{cyan!70!blue}
  \colorlet{supportABC}{violet!85!black}

  \tikzset{
    support boundary/.style={line cap=round,line join=round},
    entrance point/.style={circle,inner sep=0pt}
  }

  \draw[gray!55, line width=1.5pt] (0,0) rectangle (6,6);

  \foreach \i/\j in {
    0/0,0/1,0/2,0/4,0/5,
    1/0,1/1,1/2,1/3,1/4,1/5,
    2/0,2/1,2/2,2/3,2/4,2/5,
    3/0,3/1,3/2,3/4,
    4/0,4/1,4/2,4/3,4/4,4/5,
    5/0,5/1,5/2,5/4,5/5
  } {
    \fill[black!8] (\i,\j) rectangle ++(1,1);
  }
 \foreach \t in {1,...,7} {
    \draw[black!16, line width=.4pt] (\t,0) -- (\t,7);
    \draw[black!16, line width=.4pt] (0,\t) -- (7,\t);
  }
  \draw[support boundary, supportA, line width=3pt]
    (0,7) -- (0,1) -- (1,1) -- (1,0) -- (7,0);
  \draw[black, line width=.85pt] (0,7) -- (0,0) -- (7,0);
  \foreach \t in {0,...,7} {
    \draw[black, line width=.65pt] (\t,0) -- (\t,-.08);
    \draw[black, line width=.65pt] (0,\t) -- (-.08,\t);
  }

  \node[anchor=north, font=\small] at (0,-0.13) {$0$};
  \node[anchor=north, font=\small] at (1,-0.13) {$1$};
  \node[anchor=north, font=\small] at (2,-0.13) {$2$};
  \node[anchor=north, font=\small] at (3,-0.13) {$3$};
  \node[anchor=north, font=\small] at (4,-0.13) {$4$};
  \node[anchor=north, font=\small] at (5,-0.13) {$5$};
  \node[anchor=north, font=\small] at (6,-0.13) {$6$};
  \node[anchor=north, font=\small] at (7,-0.13) {$7$};
  \node[anchor=east, font=\small] at (-0.13,0) {$0$};
  \node[anchor=east, font=\small] at (-0.13,1) {$1$};
  \node[anchor=east, font=\small] at (-0.13,2) {$2$};
  \node[anchor=east, font=\small] at (-0.13,3) {$3$};
  \node[anchor=east, font=\small] at (-0.13,4) {$4$};
  \node[anchor=east, font=\small] at (-0.13,5) {$5$};
  \node[anchor=east, font=\small] at (-0.13,6) {$6$};
  \node[anchor=east, font=\small] at (-0.13,7) {$7$};

  \draw[support boundary, supportBCA, line width=2.2pt]
    (1,7) -- (1,2) -- (2,2) -- (2,1) -- (7,1);

  \draw[support boundary, supportBC, line width=2.2pt]
    (4,7) -- (4,2) -- (7,2);
  \draw[support boundary, supportAC, line width=2.2pt]
    (1,7) -- (1,4) -- (7,4);
  \draw[support boundary, supportABC, line width=2.2pt]
    (4,7) -- (4,6) -- (5,6) -- (5,5) -- (6,5) -- (6,4) -- (7,4);

  \foreach \p in {(0,1),(1,0)} {
    \node[entrance point, fill=supportA, minimum size=6pt] at \p {};
  }

  \foreach \p in {(1,2),(2,1)} {
    \node[entrance point, fill=supportBCA, minimum size=6pt] at \p {};
  }

  \node[entrance point, fill=supportBC, minimum size=6pt] at (4,2) {};
  \node[entrance point, fill=supportAC, minimum size=6pt] at (1,4) {};
  \foreach \p in {(4,6),(5,5),(6,4)} {
    \node[entrance point, fill=supportABC, minimum size=6pt] at \p {};
  }

  \node[anchor=south west, font=\small] at (0.10,1.10) {$[a]$};
  \node[anchor=south west, font=\small] at (1.10,2.10) {$[b],[c],[a,b]$};
  \node[anchor=south west, font=\small] at (4.10,2.10) {$[a,c]$};
  \node[anchor=south west, font=\small] at (1.10,4.10) {$[b,c]$};
  \node[anchor=south west, font=\small] at (5.10,5.10) {$[a,b,c]$};
\end{scope}

%% file: figs/tixzfigures/bifiltration_grid_sketch.tex
\providecommand{\bifiltrationfigurepath}{figs/tixzfigures}

\begin{scope}
\colorlet{ksixHighlight}{blue!12}
\colorlet{ksevenHighlight}{orange!13}
\colorlet{keightHighlight}{violet!11}
\colorlet{knineHighlight}{red!11}
\newcommand{\sketchcomplexcell}[2]{%
  \tikz[baseline=(current bounding box.center)]{%
    \fill[#1] (0,0) rectangle (1.5,1.5);%
    \input{\bifiltrationfigurepath/bifiltration_#2}%
  }%
}

\matrix (S) [
  matrix of nodes,
  nodes in empty cells,
  row sep=0pt,
  column sep=0pt,
  nodes={inner sep=0pt, outer sep=0pt}
] {
  %
  \tikz[baseline=(current bounding box.center)]{%
    \fill[ksevenHighlight] (0,0) rectangle (1.5,1.5);%
    \input{\bifiltrationfigurepath/bifiltration_K7}%
  }%
 & %
  \tikz[baseline=(current bounding box.center)]{%
    \fill[keightHighlight] (0,0) rectangle (1.5,1.5);%
    \input{\bifiltrationfigurepath/bifiltration_K8}%
  }%
 & %
  \tikz[baseline=(current bounding box.center)]{%
    \fill[keightHighlight] (0,0) rectangle (1.5,1.5);%
    \input{\bifiltrationfigurepath/bifiltration_K8}%
  }%
 & %
  \tikz[baseline=(current bounding box.center)]{%
    \fill[keightHighlight] (0,0) rectangle (1.5,1.5);%
    \input{\bifiltrationfigurepath/bifiltration_K8}%
  }%
 \\
  %
  \tikz[baseline=(current bounding box.center)]{%
    \fill[ksevenHighlight] (0,0) rectangle (1.5,1.5);%
    \input{\bifiltrationfigurepath/bifiltration_K7}%
  }%
 & %
  \tikz[baseline=(current bounding box.center)]{%
    \fill[keightHighlight] (0,0) rectangle (1.5,1.5);%
    \input{\bifiltrationfigurepath/bifiltration_K8}%
  }%
 & %
  \tikz[baseline=(current bounding box.center)]{%
    \fill[keightHighlight] (0,0) rectangle (1.5,1.5);%
    \input{\bifiltrationfigurepath/bifiltration_K8}%
  }%
 & %
  \tikz[baseline=(current bounding box.center)]{%
    \fill[keightHighlight] (0,0) rectangle (1.5,1.5);%
    \input{\bifiltrationfigurepath/bifiltration_K8}%
  }%
 \\
  %
  \tikz[baseline=(current bounding box.center)]{%
    \fill[ksixHighlight] (0,0) rectangle (1.5,1.5);%
    \input{\bifiltrationfigurepath/bifiltration_K6}%
  }%
 & %
  \tikz[baseline=(current bounding box.center)]{%
    \fill[knineHighlight] (0,0) rectangle (1.5,1.5);%
    \input{\bifiltrationfigurepath/bifiltration_K9}%
  }%
 & %
  \tikz[baseline=(current bounding box.center)]{%
    \fill[keightHighlight] (0,0) rectangle (1.5,1.5);%
    \input{\bifiltrationfigurepath/bifiltration_K8}%
  }%
 & %
  \tikz[baseline=(current bounding box.center)]{%
    \fill[keightHighlight] (0,0) rectangle (1.5,1.5);%
    \input{\bifiltrationfigurepath/bifiltration_K8}%
  }%
 \\
  %
  \tikz[baseline=(current bounding box.center)]{%
    \fill[ksixHighlight] (0,0) rectangle (1.5,1.5);%
    \input{\bifiltrationfigurepath/bifiltration_K6}%
  }%
 & %
  \tikz[baseline=(current bounding box.center)]{%
    \fill[ksixHighlight] (0,0) rectangle (1.5,1.5);%
    \input{\bifiltrationfigurepath/bifiltration_K6}%
  }%
 & %
  \tikz[baseline=(current bounding box.center)]{%
    \fill[ksevenHighlight] (0,0) rectangle (1.5,1.5);%
    \input{\bifiltrationfigurepath/bifiltration_K7}%
  }%
 & %
  \tikz[baseline=(current bounding box.center)]{%
    \fill[ksevenHighlight] (0,0) rectangle (1.5,1.5);%
    \input{\bifiltrationfigurepath/bifiltration_K7}%
  }%
 \\
};

\foreach \row in {1,...,4} {
  \foreach \col in {1,...,4} {
    \draw[black!70, line width=0.65pt]
      (S-\row-\col.north west) rectangle (S-\row-\col.south east);
  }
}


\end{scope}

%% file: figs/tixzfigures/directed_parameter_space_complex.tex
\begin{scope}[x=1.35cm,y=1.35cm]
  \colorlet{supportBC}{cyan!70!blue}
  \colorlet{supportAC}{red!80!black}
  \colorlet{pathAlpha}{green!60!black}
  \colorlet{pathBeta}{magenta!75!black}

  \tikzset{
    support boundary/.style={line cap=round,line join=round},
    entrance point/.style={circle,inner sep=0pt},
directed path/.style={
  line width=2.35pt,
  line cap=round,
  line join=round,
  opacity=.45
}
  }

  \foreach \i/\j in {
    0/0,0/2,0/3,
    1/1,1/2,1/3,
    2/0,2/1,2/2,2/3,
    3/0,3/1,3/2,3/3
  } {
    \fill[black!8] (\i,\j) rectangle ++(1,1);
  }

  \foreach \t in {1,...,4} {
    \draw[black!16, line width=.4pt] (\t,0) -- (\t,4);
    \draw[black!16, line width=.4pt] (0,\t) -- (4,\t);
  }

  \draw[gray!55, line width=1.5pt] (0,0) rectangle (4,4);
  \node[anchor=south west, font=\small] at (0.10,4.08) {$C$};

  \draw[support boundary, supportBC, line width=1.6pt]
    (0,4) -- (0,2) -- (2,2) -- (2,0) -- (4,0);
  \draw[support boundary, supportAC, line width=1.6pt]
    (1,4) -- (1,1) -- (4,1);

  \draw[black, line width=.85pt] (0,4) -- (0,0) -- (4,0);
  \foreach \t in {0,...,4} {
    \draw[black, line width=.65pt] (\t,0) -- (\t,-.08);
    \draw[black, line width=.65pt] (0,\t) -- (-.08,\t);
  }

  \foreach \t in {0,...,4} {
    \node[anchor=north, font=\small] at (\t,-0.13) {$\t$};
    \node[anchor=east, font=\small] at (-0.13,\t) {$\t$};
  }

  \draw[directed path, pathAlpha] (0,0) -- (0,2) -- (4,2) -- (4,4);
  \node[anchor=south, font=\small, text=pathAlpha] at (2.15,2.08) {$\alpha$};

  \draw[directed path, pathBeta] (0,0) -- (2,0) -- (2,4) -- (4,4);
  \node[anchor=west, font=\small, text=pathBeta] at (2.08,3.12) {$\beta$};

  \node[entrance point, fill=black!75, minimum size=4.5pt] at (0,0) {};
  \node[anchor=south west, font=\small] at (0.10,0.10) {$[a],[b],[c]$};

  \foreach \p in {(0,2),(2,0)} {
    \node[entrance point, fill=supportBC, minimum size=5pt] at \p {};
  }
  \node[entrance point, fill=supportAC, minimum size=5pt] at (1,1) {};
  \node[anchor=south west, font=\small] at (2.10,0.10) {$[b,c]$};
    \node[anchor=south west, font=\small] at (0.10,2.10) {$[b,c]$};
  \node[anchor=south west, font=\small] at (1.10,1.10) {$[a,c]$};
\end{scope}

%% file: figs/tixzfigures/directed_alpha_beta_filtrations.tex
\providecommand{\bifiltrationfigurepath}{figs/tixzfigures}

\begin{scope}
\colorlet{ksixHighlight}{blue!12}
\colorlet{ksevenHighlight}{orange!13}
\colorlet{keightHighlight}{violet!11}

\newcommand{\directedfiltrationcell}[2]{%
  \tikz[baseline=(current bounding box.center)]{%
    \fill[#1] (0,0) rectangle (1.5,1.5);%
    \input{\bifiltrationfigurepath/bifiltration_#2}%
  }%
}

\matrix (F) [
  matrix of nodes,
  nodes in empty cells,
  row sep=0pt,
  column sep=0pt,
  nodes={inner sep=0pt, outer sep=0pt}
] {
  \tikz[baseline=(current bounding box.center)]{%
    \fill[ksixHighlight] (0,0) rectangle (1.5,1.5);%
    \input{\bifiltrationfigurepath/bifiltration_K6}%
  }%
 &
  \tikz[baseline=(current bounding box.center)]{%
    \fill[ksevenHighlight] (0,0) rectangle (1.5,1.5);%
    \input{\bifiltrationfigurepath/bifiltration_K7}%
  }%
 &
  \tikz[baseline=(current bounding box.center)]{%
    \fill[keightHighlight] (0,0) rectangle (1.5,1.5);%
    \input{\bifiltrationfigurepath/bifiltration_K8}%
  }%
 \\
};

\foreach \col in {1,...,3} {
  \draw[black!70, line width=0.65pt]
    (F-1-\col.north west) rectangle (F-1-\col.south east);
}
\end{scope}